\documentclass[11pt]{article}
\usepackage{amssymb}
\parskip=\medskipamount
\font\fr=eufm10  scaled \magstep 1   
\font\ddpp=msbm10  scaled \magstep 1  

\def\QED{\hskip0.1em\hfill\null\ \null\nobreak\hfill
\kern3pt\lower1.8pt\vbox{\hrule\hbox   {\vrule\kern1pt\vbox{\kern1.7pt
\hbox{$\scriptstyle   QED$}\kern0.2pt}\kern1pt\vrule}\hrule}}
\def\R{\hbox{\ddpp R}}               
\def\S{\hbox{\ddpp S}}
\def\L{\hbox{\ddpp L}}

\begin{document}
\title{{\bf On the geometry of generalized Chaplygin systems}}     
\author{
{\sc By FRANS CANTRIJN}\\
{\em Department of Mathematical Physics and Astronomy, Ghent University}\\
{\em Krijgslaan 281, B-9000 Ghent, Belgium}\\
{\em e-mail}: Frans.Cantrijn@rug.ac.be\\
{\sc JORGE CORT\'ES} \dag \qquad {\sc MANUEL DE LE\'ON} \ddag\\
{\sc and DAVID MART\'IN DE DIEGO} \P \\
{\em Laboratory of Dynamical Systems, Mechanics and Control} \\
{\em Instituto de Matem\'aticas y F{\'\i}sica Fundamental, CSIC} \\
{\em Serrano 123, 28006 Madrid, Spain}\\
\dag {\em e-mail}: j.cortes@imaff.cfmac.csic.es\quad 
\ddag {\em e-mail}: mdeleon@imaff.cfmac.csic.es \\ 
\P {\em e-mail}: d.martin@imaff.cfmac.csic.es}
\date{\today}
\maketitle     

Short title: {\em Chaplygin systems}\\
2000 MS Classification: 37J60, 70Hxx

\bigskip
\newpage
\begin{center}
{\it Abstract}
\end{center}
Some aspects of the geometry and the dynamics of generalized Chaplygin systems 
are investigated. First, two different but complementary approaches to the 
construction of the reduced dynamics are reviewed: a symplectic approach and 
an approach based on the theory of affine connections. Both are mutually compared 
and further completed. Next, a necessary and sufficient condition is derived for 
the existence of an invariant measure for the reduced dynamics of generalized 
Chaplygin systems of mechanical type. A simple example is then constructed of 
a generalized Chaplygin system which does not verify this condition, thereby 
answering in the negative a question raised by Koiller.

\begin{center}
1. {\it Introduction}
\end{center}
In the past two decades, the rich developments in the field of 
Geometric Mechanics have led to a considerable progress 
in the study of the geometrical structure and of the dynamics of 
mechanical systems with nonholonomic constraints. Several concepts and 
techniques, familiar from the geometric approach to Lagrangian and Hamiltonian
mechanics, have been succesfully adapted or extended to the framework of 
nonholonomic mechanics. Examples of this succesful ``transfer of ideas" 
can be found, among others, in the symplectic and Poisson descriptions of 
constrained systems \cite{BaSn,CuKeSnBa,KoMa1,KoMa3,LeMaMa,LeMa,VaMa} and, in
particular, in the study of nonholonomic systems with symmetry (reduction and 
reconstruction of the dynamics, stability of relative equilibria ...) where 
elements are being used from the theory of symplectic reduction and from the 
theory of principal connections and Ehresmann connections \cite{BaSn,BlKrMaMu,CaLeMaMa,CoLe,Ko,Ma}.  
Whereas most of the foregoing treatments are concerned with the autonomous case
(i.e.\ time-invariant systems with time-independent constraints), the extension
to time-dependent nonholonomic systems, using the geometric formalism of
jet bundle theory, has also been treated by several authors: see 
e.g.\ \cite{Krup,Massa,Sau} and references therein.

An important topic which is receiving growing attention in the literature, 
concerns the identification and characterization of a suitable notion of 
complete integrability of nonholonomic systems (see e.g.\ 
\cite{Ar,BaCu,BlCr,DrGaJo,Jo,Koz,VeVe}). As is well known, an (unconstrained) 
Hamiltonian system on a $2n$-dimensional phase space is called completely 
integrable if it admits $n$ independent integrals of motion in involution. 
It then follows from the Arnol'd-Liouville theorem that, when assuming compactness 
of the common level sets of these first integrals, the motion in phase space 
is quasi-periodic and consists of a winding on $n$-dimensional invariant tori 
(see e.g.\ \cite{Ar}, Chapter 4).
For the integrability of a nonholonomic system with $k$ 
constraints one needs, in general, $2n-k-1$ independent first integrals. 
It turns out, however, that for a nonholonomic system
which admits an invariant measure, ``only" $2n-k-2$ first 
integrals are needed in order to reduce its integration to quadratures, 
and in such a case - again assuming compactness of the common level sets
of the first integrals - the phase space trajectories of the system live 
on 2-dimensional invariant tori \cite{Ar}.
Several authors have studied the problem of the existence of invariant 
measures for some special classes of nonholonomic systems. For instance, 
Veselov and Veselova \cite{VeVe} have studied nonholonomic geodesic flows 
on Lie groups with a left-invariant metric and a right-invariant nonholonomic 
distribution (the so-called LR systems). Kozlov \cite{Koz} has treated 
the analogous problem for left-invariant constraints. Their results have 
been very useful for finding new examples of completely integrable 
nonholonomic dynamical systems \cite{DrGaJo,Jo,VeVe}.

In the present paper we will mainly be concerned with a particular, but 
important class of nonholonomic systems: the so-called generalized 
Chaplygin systems. A system of generalized Chaplygin type is described by 
a principal connection on a principal fiber bundle $Q \longrightarrow Q/G$, 
whose horizontal distribution determines the constraint submanifold, and 
a $G$-invariant Lagrangian $L:TQ \longrightarrow \R$. Generalized Chaplygin 
systems exhibit a very nice geometric structure, which has been discussed 
extensively in the recent literature 
(see e.g.\ \cite{BlKrMaMu,CaLeMaMa2,CoLe,Ko,LeMa,St,Ya}, and references therein). 
One of the peculiarities of these systems is that, after reduction, they take 
on the form of an unconstrained system, subject to an ``external" force 
of a special type. In \cite{St} it is shown that in the abelian Chaplygin case, 
the reduced equations can be rewritten in the form of a Hamiltonian system with 
respect to an almost symplectic structure. We will see that this also holds
in the nonabelian case. 

In his pioneering paper on the reduction of nonholonomic systems 
with symmetry, Koiller touches the problem of the existence of an
invariant measure for the reduced dynamics of generalized Chaplygin systems 
(see \cite{Ko}, Section 9). Based on several known examples of such systems 
which do admit an invariant measure, Koiller asks the question whether 
this property might perhaps hold in general. The main result 
of the present paper is the derivation of a necessary and sufficient 
condition for the existence of an invariant measure 
for the reduced dynamics of a generalized Chaplygin system whose Lagrangian 
is of pure kinetic energy type. This condition then enables us to 
give a negative answer to Koiller's question by constructing a simple counter example.

Clearly, there is still much work to be done in order to obtain a deeper
and more complete understanding of the structure and behaviour of completely 
integrable nonholonomic systems. Any further progress into this matter will 
also be of importance for the study of perturbations of integrable nonholonomic 
systems. Concerning the particular case of generalized Chaplygin systems, it 
would be of interest, for instance, to identify general classes of systems which 
do admit an invariant measure, and to use their characterization to tackle the 
problem of complete integrability.

This paper is organized as follows. In Section 2, we first give a very brief 
review of some aspects of the theory of connections on principal bundles  
before introducing the notion of generalized Chaplygin system. Sections 3 and 4 
are devoted to two alternative approaches to the description of generalized 
Chaplygin systems, with special emphasis on the reduction and reconstruction 
of their dynamics: a symplectic approach and an approach based on 
the use of affine connections. In Section 5, two examples are treated which 
exhibit a different behaviour in the reduction process. Section 6 then deals with 
the relation between the symplectic and the affine connection approach. The existence
of an invariant measure for the reduced equations of generalized Chaplygin systems 
of mechanical type is studied in Section 7. The main result 
of the paper (Theorem 7.5) presents a necessary and sufficient condition for the 
existence of such an invariant measure, and an example is given of a system for 
which this condition is not satisfied. Finally, some particular results 
concerning the theory of affine connections are discussed in an Appendix.

Throughout this paper, we are working in the category of $C^{\infty}$-manifolds 
(with smooth maps, tensor fields, etc...). For convenience, we shall usually not 
make a notational distinction between a vector bundle over a manifold and
the module of its smooth sections, i.e. if $\pi:F \rightarrow N$ denotes a 
vector bundle over a manifold $N$ (for instance a subbundle of $TN$), then 
$X \in F$ simply means that $X$ is a smooth section of $\pi$.
The sole exception will be the sporadic use of the notation ${\frak X}(N)$
for the module of vector fields on $N$.
If $\cal D$ is distribution on a manifold $N$, then its annihilator is a 
codistribution on $N$ denoted by ${\cal D}^o$. Both $\cal D$ and ${\cal D}^o$ 
will also be identified with the corresponding vector subbundles of $TN$, resp. 
$T^*N$. Finally, the tangent map of a mapping $f$ between manifolds will be 
denoted interchangeably by $Tf$ and $f_*$.
 
\begin{center} 
2. {\it Generalized Chaplygin systems}
\end{center}

In this section, we first briefly review the notion of a principal
connection on a principal fibre bundle. For details we refer to \cite{KoNo}. 
Secondly, we consider a geometric framework for general nonholonomic 
mechanical systems and then introduce the concept of generalized Chaplygin 
systems. 

2.1. {\it Principal connections}

Let $Q$ be the configuration manifold of a physical system and assume that
there is a left action of a Lie group $G$ on $Q$: 
\[
\begin{array}{rrcl}
\Psi:& G\times Q&\longrightarrow&Q\\
       & (g, q)&\longmapsto&\Psi(g,q)=\Psi_g(q)=gq \, .
\end{array}
\]
Note that we consider here left actions, which is the usual convention in
mechanics, instead of the right ones, as considered for instance in 
\cite{KoNo}. The orbit through a point $q$ is 
$\hbox{Orb}_G(q)=\{gq\; |\; g\in G\}$. We denote by $\hbox{\fr g}$ 
the Lie algebra of $G$. For any element $\xi\in \hbox{\fr g}$, 
$\xi_Q$ will denote the corresponding infinitesimal generator of the group
action on $Q$. Then, 
\[
T_q(\hbox{Orb}_G(q))=\{\xi_Q (q)\; |\; \xi\in \hbox{\fr g}\} \, .
\]
Assuming that the action $\Psi$ is free and proper, we can endow the
quotient space $Q/G=M$ with a manifold structure such that the canonical
projection $\pi: Q \longrightarrow M$ is a surjective submersion. In the
framework of the mechanics of (coupled) rigid bodies, for instance,
the quotient manifold $M$ is commonly called the {\bf shape space} of the
system under consideration. We then have that $Q(M, G, \pi)$ is a principal 
bundle with bundle space $Q$, base space $M$, structure group $G$ and 
projection $\pi$. Note that the kernel of $\pi_* (= T\pi)$ consists of the 
vertical tangent vectors, i.e.\ the vectors tangent to the orbits of $G$ 
in $Q$. We shall denote the bundle of vertical vectors by ${\cal V}_{\pi}$, with 
$({\cal V}_{\pi})_q=T_q (\hbox{Orb}_G(q))$, $q \in Q$.

A principal connection on $Q(M, G, \pi)$ can be defined as 
a distribution ${\cal H}$ on $Q$ satisfying: 
\begin{enumerate}
\item $T_qQ= {\cal H}_q \oplus ({\cal V}_{\pi})_q$, $\forall q \in Q$;
\item ${\cal H}_{gq}=T_q\Psi_g({\cal H}_q)$, i.e. the distribution
${\cal H}$ is $G$-invariant; 
\item ${\cal H}_q$ depends smoothly on $q$.
\end{enumerate}
The subspace ${\cal H}_{q}$ of $T_qQ$ is called the horizontal subspace at $q$ 
determined by the connection. Alternatively, a principal connection can 
be characterized by a $\hbox{\fr g}$-valued 1-form $\gamma$ on $Q$ 
satisfying the following conditions 
\begin{enumerate}
\item $\gamma(\xi_Q(q)) = \xi$ for all $\xi \in \hbox{\fr g}$,
\item $\gamma(T\Psi_g X) = \hbox{Ad}_g (\gamma(X))$ for all $X \in TQ$.
\end{enumerate}
The horizontal subspace at $q$ is then given by ${\cal H}_{q}=\{v_q \in
T_q Q \; | \; \gamma(v_q)=0 \}$. A vector field $X$ on $Q$ is called horizontal
if $X(q) \in {\cal H}_q$ at each point $q$. 

Given a principal connection, we have that every vector $v \in T_q Q$ can
be uniquely written as 
\[
v=v_1+ v_2\, ,
\]
with $v_1\in {\cal H}_q$ and $v_2\in ({\cal V}_{\pi})_q$. We denote by ${\bf h}:TQ
\longrightarrow {\cal H}$ and ${\bf v}:TQ \longrightarrow {\cal V}_{\pi}$ the
corresponding horizontal and vertical projector, respectively.
The horizontal lift of a vector field $Y$ on $M$ is the unique vector
field $Y^h$ on $Q$ which is horizontal and projects onto $Y$. 

The curvature $\Omega$ of the principal connection is the $\hbox{\fr g}$-valued
2-form on $Q$ defined as follows: for each $q \in Q$ and $u,v \in T_qQ$ 
\[
\Omega(u,v) = d\gamma({\bf h}u,{\bf h}v) = 
- \gamma([U^h,V^h]_q) \, ,
\]
where $U^h$ and $V^h$ are the horizontal lifts of any two (local) vector fields 
$U$ and $V$ on $M$ for which $U^h(q) = {\bf h}u$ and $V^h(q) = {\bf h}v$,
respectively. The curvature measures the lack of 
integrability of the horizontal distribution and plays a fundamental role in 
the theory of holonomy (see \cite{KoNo} for a comprehensive treatment). 

2.2. {\it Nonholonomic mechanics and generalized Chaplygin systems}

A nonholonomic Lagrangian system consists of a Lagrangian $L$ defined on
$TQ$, the tangent bundle of an $n$-dimensional configuration manifold
$Q$, and constraints which determine a submanifold ${\cal D}$ of $TQ$. This 
means that the only allowable velocities are those belonging to ${\cal D}$.
In case ${\cal D}$ is a vector subbundle of $TQ$, we are dealing with the case
of linear constraints and if, in addition, this subbundle corresponds to 
an integrable distribution, we are reduced to the case of holonomic constraints. 
In the sequel, we will always assume that $\tau_Q({\cal D})=Q$, where 
$\tau_Q: TQ \rightarrow Q$ is the tangent bundle projection. 

Let $(q^A)$, $A=1,...,n$, be local coordinates on $Q$ and denote the induced 
bundle coordinates on $TQ$ by $(q^A,\dot{q}^A)$. In a local description, a 
constraint submanifold ${\cal D}$ of codimension $k$ can be defined by the
vanishing of $k$ independent functions $\phi_i$ (the constraint functions). 
In the following, we shall only consider the case of linear constraints, such that
the functions $\phi_i$ can be taken to be of the form $\phi_i(q,\dot{q}) = 
\mu_{iA}(q)\dot{q}^A$, $i=1,...,k$. Application of d'Alembert's principle, using
the common notion of virtual displacement for the case of linear constraints, 
then leads to the constrained equations of motion
\begin{equation}\label{noholo}
\frac{d}{dt} \left( \frac{\partial L}{\partial \dot{q}^A} \right) 
- \frac{\partial L}{\partial q^A} = \lambda^i \mu_{iA} \, ,
\end{equation}
which, together with the constraint equations $\mu_{iA}\dot{q}^A=0$ 
($i=1,...,k$), determine the dynamics of the nonholonomic system. Here, the
$\lambda^i$ are Lagrange multipliers to be determined. The right-hand side of 
equation (\ref{noholo}) precisely represents the ``reaction force" induced by
the constraints.

We now describe the structure of a so-called {\bf generalized Chaplygin system} 
\cite{BlKrMaMu,Ko,LeMa}. The configuration manifold $Q$ of a generalized Chaplygin
system is a principal $G$-bundle $\pi : Q \longrightarrow Q/G$, and the 
constraint submanifold ${\cal D}$ is given by the horizontal distribution 
${\cal H}$ of a principal connection $\gamma$ on $\pi$. Furthermore, 
the system is described by a regular Lagrangian $L : TQ \longrightarrow \R$, 
which is $G$-invariant for the lifted action of $G$ on $TQ$. 
In this paper, we shall mainly restrict our attention to systems of mechanical 
type for which $L=T-V$, where $T:TQ \longrightarrow \R$ 
is the kinetic energy, corresponding to a Riemannian metric $g$ on $Q$, 
and $V:Q \longrightarrow \R$ is the potential energy. 
Whenever we add the word ``mechanical" to the description of a system, 
we shall always be referring to this situation. For the case of a Chaplygin system
we then suppose, in addition, that both the potential energy and the metric 
$g$ are $G$-invariant so that
\[
{\cal L}_{\xi_Q} V = 0, \; \; {\cal L}_{\xi_Q}g = 0 \, ,
\]
for all $\xi \in \hbox{\fr g}$. In particular, it follows that all fundamental 
vector fields $\xi_Q$ are Killing vector fields.

Typical problems in mechanics, such as the vertical and
the inclined rolling disk, the nonholonomic free particle and the two
wheeled carriage, can be interpreted as generalized Chaplygin systems in 
the above sense. Systems of that type also occur in many problems of 
robotic locomotion \cite{KeMu} and motions of microorganisms at low Reynolds 
number \cite{ShWi}. The dynamics of generalized Chaplygin systems can be 
described from a symplectic point of view, as we will show in the next
section. However, for generalized Chaplygin systems of mechanical type 
there also exists a nice geometric description in terms of affine connections. 
This will be outlined in Section 4. 

We conclude this section with some comments concerning terminology. 
Classically, a mechanical system with Lagrangian $L(q^A,\dot{q}^A)$, 
$A = 1, \ldots, n$, subject to $k$ linear nonholonomic constraints, is said to 
be of Chaplygin type if coordinates $(q^{a},q^{\alpha})$ can be found, 
with $a = 1, \ldots, k$ and $\alpha = k+1, \ldots, n$, such that the constraints 
can be written in the form 
$\dot{q}^a = B^a_{\alpha}(q^{k+1}, \ldots, q^{n})\dot{q}^{\alpha}$
and such that $L$ does not depend on the coordinates $q^a$ (see e.g.\ \cite{NeFu}). 
Such a system can be (locally) interpreted as a special case of the 
generalized Chaplygin systems introduced above, with $Q = \R^n$ and with an 
action defined by the abelian group $G =\R^{k}$ (cf.\ \cite{Ko}). Koiller
refers to the more general case, considered in the present paper, 
as ``non-abelian Chaplygin systems". In the literature on nonholonomic systems 
with symmetry, these systems are also said to be of ``principal" or ``purely 
kinematical" type \cite{BlKrMaMu,CaLeMaMa,Ko}. Finally, it should be emphasized 
that the generalized Chaplygin systems studied in \cite{LeMa} are still of a more 
general type than the ones we consider here in that they are defined
on fibre bundles which need not be principal bundles.

\begin{center}
3. {\it Symplectic approach}
\end{center}

As mentioned in the Introduction, considerable efforts have 
been made to adapt and extend several ideas and techniques from the geometric 
treatment of unconstrained problems to the study of systems with nonholonomic 
constraints. The subject has been approached from different points of view: 
a Lagrangian approach \cite{BlKrMaMu,KoMa1,Os}, a Hamiltonian approach 
\cite{BaSn,CuKeSnBa} and a formulation in terms of (almost-)Poisson structures \cite{Ca,KoMa3,Ma,VaMa}. 
In this section, we start with a brief review of an interesting symplectic 
approach to nonholonomic dynamics, developed in \cite{LeMaMa,LeMa} which, 
in particular, is well suited for the treatment of nonholonomic systems with 
symmetry. We first outline this approach for general nonholonomic systems, and 
then turn to the case of generalized Chaplygin systems for the 
discussion of reduction and reconstruction of the constrained dynamics.

We start by fixing some notations. In terms of the tangent bundle coordinates
$(q^A,\dot{q}^A)$, let us denote by 
$\displaystyle{\Delta = \dot{q}^A \frac{\partial} {\partial \dot{q}^A}}$ 
the dilation vector field on $TQ$ and by
$\displaystyle{S = dq^A \otimes \frac{\partial}{\partial \dot{q}^A}}$
the canonical vertical endomorphism (see \cite{LeRo}). The action of $S$
on a 1-form will be denoted by $S^*$. Then we can
define the Poincar\'e-Cartan 1-form and 2-form, corresponding to a given
Lagrangian $L$, by $\theta_L=S^*dL$ and $\omega_L = -d\theta_L$, respectively. 
We further have that $E_L = \Delta L-L$
represents the energy function of the system. If the Lagrangian $L$ is regular,
which will always be tacitly assumed in the sequel, $\omega_L$ is symplectic 
and induces two isomorphisms of $C^{\infty}(TQ)$-modules: 
\[ \flat_L : \hbox{\fr X}(TQ) \longrightarrow \Omega^1(TQ) \, , \quad \sharp_L :
\Omega^1(TQ) \longrightarrow \hbox{\fr X}(TQ) \, , 
\] 
where $\flat_L(X) = i_X \omega_L$ and $\sharp_L = \flat_L^{-1}$. In the absence 
of constraints, the dynamics of the Lagrangian system, with Lagrangian $L$, 
is given by the (unique) solution $\Gamma_L$ of the
equation $i_{\Gamma_L}\omega_L = dE_L$, i.e. $\Gamma_L = \sharp_L(dE_L)$. 
Indeed, $\Gamma_L$ is a second-order differential
equation field (SODE for short) whose integral curves $(q^A(t),\dot{q}^A(t) \equiv
\frac{dq^A}{dt}(t))$ are determined by the solutions $q^A(t)$ of
the Euler-Lagrange equations for $L$. 

In the presence of nonholonomic constraints, the equations of motion must be 
modified in order to incorporate the constraints into the picture. 
Since we confine ourselves here to the case of linear constraints, the constraint 
submanifold ${\cal D}$ is a vector subbundle of $TQ$, determined by a 
(nonintegrable) distribution on $Q$ which we also denote by ${\cal D}$.
In addition, we shall always assume that the constraints verify the so-called 
``admissibility condition" (see e.g.\ \cite{LeMa}), i.e.\ for all $x \in {\cal D}$
\[
\hbox{dim}\, (T_x{\cal D})^o = \hbox{dim} \, S^* ((T_x{\cal D})^o) \, ,
\]
where the annihilator of $T_x{\cal D}$ is taken in $T^*_xTQ$. Locally, ${\cal D}$
is described by equations of the form
\[
\mu_{iA}(q)\dot{q}^A = 0, \qquad i = 1, \ldots k,
\]
with $k = \hbox{codim}\,({\cal D})$. 
Next, we define a distribution $F$ on $TQ$ along ${\cal D}$, by prescribing 
its annihilator to be a subbundle of $T^*TQ_{|{\cal D}}$ which represents the 
bundle of reaction forces. More precisely, we set $F^{o} = S^*((T{\cal D})^{o})$. 
The equations of motion for the nonholonomic system are then given by
\begin{equation}
\label{equat3}
\left\{
\begin{array}{rcl}
&& (i_{X}\omega_{L} - dE_{L})_{|{\cal D}} \in F^{o} \; ,\\
&& X_{|{\cal D}} \in T{\cal D} \; .
\end{array}\right.
\end{equation}
This system will have a unique solution $X$ provided the ``compatibility 
condition" holds, i.e.\ $F^\perp \cap T{\cal D} =0$, where 
$F^\perp = \sharp_L (F^o)$. For systems of mechanical type, with linear 
nonholonomic constraints and positive definite kinetic energy, this condition 
is always fulfilled. Indeed, putting 
\[
C_{ij}= -\mu_{iA}W^{AB}\mu_{jB}, \quad i,j =1, \ldots k,
\]
where $(W^{AB})$ is the inverse of the Hessian matrix 
$\displaystyle{\left( \frac{\partial^2 L}{\partial  \dot{q}^A \partial 
\dot{q}^B} \right)}$, compatibility locally translates into 
regularity of the matrix $(C_{ij})$ (cf.\ \cite{LeMa}, where the 
compatibility condition for a nonholonomic system was called the 
regularity condition). Under this condition, one can show that we 
have a direct sum decomposition of $T_{\cal D}TQ = T{\cal D} \oplus F^{\perp}$
and that the constrained dynamics $X$ is obtained by
projecting the unconstrained Euler-Lagrange vector field $\Gamma_L$ (restricted
to ${\cal D}$) onto $T{\cal D}$ with respect to this decomposition.

It should be pointed out that the solution $X$ of (\ref{equat3}) satisfies
automatically the SODE condition along ${\cal D}$, i.e.\ $S(X)_{|{\cal D}} =
\Delta_{|{\cal D}}$. This implies that, in local coordinates, the integral 
curves of $X$ on ${\cal D}$ are of the form $(q^A(t), \dot{q}^A(t) \equiv
\frac{dq^A}{dt}(t))$, whereby the $q^A(t)$ are solutions of the system of 
differential equations (\ref{noholo}), together with the constraint equations 
$\mu_{iA}(q) \dot{q}^A = 0$, $i=1,...,k$. The local coordinate expression for 
$X$ reads
\begin{eqnarray*}
X &=& \dot{q}^A \frac{\partial}{\partial q^A} \\
&+& W^{AB}\left( \frac{\partial L}{\partial q^B}- 
\frac{\partial p_B}{\partial q^C}\dot{q}^C +
C^{ij}\frac{\partial \mu_{iD}}{\partial q^C}\mu_{jB}\dot{q}^C\dot{q}^D +  
W^{CD}\left( \frac{\partial L}{\partial q^D}-\dot{q}^E\frac{\partial p_D}
{\partial q^E} \right)C^{ij}\mu_{jB} \mu_{iC} \right) 
\frac{\partial}{\partial \dot{q}^A} \, ,
\end{eqnarray*}
where $(C^{ij})$ is the inverse of the matrix $(C_{ij})$ introduced above and where, 
for ease of writing, we have put 
$p_A = \displaystyle{\frac{\partial L}{\partial \dot{q}^A}}$, $A=1,...,n$.

Before proceeding, we recall that for nonholonomic Lagrangian systems with 
constraints which are linear (or, more general, homogeneous) in the velocities,
the energy $E_L$ is a conserved quantity (see e.g.\ \cite{CaLeMaMa}). This, 
therefore, in particular applies to the generalized Chaplygin systems 
considered in this paper.

3.1. {\it Reduction} 

The reduction theory of nonholonomic systems with symmetry, and 
related aspects, has become an intensive field of research \cite{BaSn,BlKrMaMu,CaLeMaMa,CaLeMaMa2,CoLe,KoMa3,Ma,Os}.
Here we consider the special case of nonholonomic systems of 
generalized Chaplygin type.  

The given data are (cf.\ Section 2): a principal $G$-bundle 
$\pi: Q \longrightarrow M = Q/G$, associated to a free and proper action
$\Psi$ of $G$ on $Q$, a Lagrangian $L: TQ \longrightarrow \R$ which is $G$-invariant
with respect to the lifted action on $TQ$, and linear nonholonomic constraints 
determined by the horizontal distribution (here denoted as ${\cal D}$) of a 
principal connection $\gamma$ on $\pi$. Taking into account the available 
symmetries, we can reduce the number of degrees of freedom of the problem. 
In the following, we review the various geometric concepts involved in the 
symplectic approach to this reduction process \cite{CaLeMaMa,CaLeMaMa2,CoLe}. 

{\bf The symplectic action.} Consider the lifted action of $G$ on $TQ$, i.e.\
$\hat{\Psi}: G \times TQ \longrightarrow TQ$ with $\hat{\Psi}(g,v_q) = 
T\Psi_g (v_q) (= {\Psi_g}_*(v_q))$ for any $g \in G$ and $v_q \in T_qQ$. This 
action is also free and proper and, moreover, it is symplectic with
respect to $\omega_L$. For any $\xi \in \hbox{\fr g}$, the infinitesimal generators 
$\xi_{TQ}$ and $\xi_Q$ of $\hat{\Psi}$ and $\Psi$, respectively, are 
$\tau_Q$-related, i.e.  
\begin{equation}\label{xi}
{\tau_Q}_* \circ \xi_{TQ} = \xi_Q \circ \tau_Q\;.
\end{equation}
Let us denote by $\rho:TQ \longrightarrow \overline{TQ} = TQ/G$ the natural 
projection. From the given assumptions it follows that the constraint
submanifold ${\cal D}$, the energy $E_L$ and the vector subbundle 
$F$ are $G$-invariant. The induced action of $G$ on ${\cal D}$, 
i.e.\ the restriction of $\hat{\Psi}$ to $G \times {\cal D}$, is still free 
and proper and we can regard the orbit space $\bar{\cal D} = {\cal D}/G$ as 
a submanifold of $TQ/G$. Note that there exists a natural identification 
${\cal D} \cong Q \times_{Q/G} T(Q/G)$ as principal $G$-bundles over $T(Q/G)$.
The isomorphism is obtained by mapping $v_q \in {\cal D}$ onto $(q,{\pi}_*(v_q))$.
It then follows that ${\cal D}/G$ can be naturally identified with $T(Q/G)$ 
and we have, in particular,
\begin{equation}\label{proj}
\rho_{|{\cal D}} = {\pi_*}_{|{\cal D}}.
\end{equation} 
Henceforth, the restriction of $\rho$ to ${\cal D}$ will also be simply denoted 
by $\rho$. 

{\bf The connection.} A direct computation shows that along ${\cal D}$ 
we have that ${\cal V}_\rho \cap F = 0$, where ${\cal V}_\rho$ denotes the 
subbundle of $TTQ$ which is vertical with respect to the projection $\rho$, 
i.e.\ ${\cal V}_\rho=\hbox{Ker}\, T\rho$. Observing that 
${{\cal V}_{\rho}}_{|{\cal D}} \subset T{\cal D}$, it then easily follows that  
\[
T{\cal D} = (F \cap T{\cal D}) \oplus  {{\cal V}_{\rho}}_{|{\cal D}} \; .
\]
Since ${\cal U}=F \cap T{\cal D}$ is $G$-invariant, the above decomposition of
$T{\cal D}$ defines a principal connection $\Upsilon$ on the principal
$G$-bundle $\rho:{\cal D} \longrightarrow \bar{\cal D}={\cal D}/G$,
with horizontal subspace ${\cal U}_x$ at each $x\in {\cal D}$. Moreover, we
have that the dynamics $X$ of the generalized Chaplygin system, determined
by (\ref{equat3}), belongs to ${\cal U}$. We shall denote the projection of 
the bundle ${\cal U}$ under $T\rho (= \rho_*)$ by $\bar{\cal U}$.

We can summarize the situation in the following diagram:
\[
\begin{array}{ccc}
T{\cal D} = {\cal U} \oplus {\cal V}_\rho & \stackrel{\rho_*}{\longrightarrow} 
& T\bar{\cal D} \cong \bar{\cal U} \\ 
\downarrow & & \downarrow \\
{\cal D} & \stackrel{\rho} {\longrightarrow} & T(Q/G) \cong \bar{\cal D}  
\end{array}
\]

The principal connection $\Upsilon$ is obviously related to the
original connection $\gamma$ of the Chaplygin system. Indeed, take
$w \in T_{v_q}{\cal D}$ and consider ${\tau_Q}_*w \in T_qQ$. Then, we
can write 
\[
{\tau_Q}_*w = \left({\tau_Q}_*w - (\gamma({\tau_Q}_*w))_Q(q)\right) + 
(\gamma({\tau_Q}_*w))_Q(q) \, ,
\]
where ${\tau_Q}_*w - (\gamma({\tau_Q}_*w))_Q(q) \in {\cal D}_q$ and
$(\gamma({\tau_Q}_*w))_Q(q) \in ({\cal V}_{\pi})_q$. Putting 
$\gamma({\tau_Q}_*w) = \xi \in \hbox{\fr g}$, a direct computation shows 
that $w - \xi_{TQ}(v_q) \in {\cal U}$ and, consequently, 
$w= (w-\xi_{TQ}(v_q)) +\xi_{TQ}(v_q)$ is the $({\cal U},{\cal V}_\rho)$ 
decomposition of $w$. Herewith we have proved the following property.

{\sc Proposition 3.1.} {\it The connection 1-forms $\Upsilon$ and $\gamma$ are 
related by $\Upsilon = \tau_Q^*\gamma$, i.e.\ 
$\Upsilon_{v_q}(w) = \gamma_q({\tau_Q}_*w)$ for any $v_q \in {\cal D}$ and 
$w \in T_{v_q}{\cal D}$.}

Let us denote the horizontal projectors, associated to $\gamma$, resp.\ $\Upsilon$, 
by ${\bf h_{\gamma}} : TQ \longrightarrow {\cal D} (\subset TQ)$, resp.\ 
${\bf h_{\Upsilon}}: T{\cal D} \longrightarrow {\cal U} (\subset T{\cal D})$. 
Likewise, the vertical projectors onto ${\cal V}_{\pi}$, resp. ${\cal V}_{\rho}$,
will be denoted by ${\bf v_{\gamma}}$, resp. ${\bf v_{\Upsilon}}$. In order
not to further overload the notations, we will use the same superscript ${}^h$
for the horizontal lifts of vectors (vector fields) with
respect to either $\gamma$ or $\Upsilon$; in principle it should always be 
clear from the context which horizontal lift operation is being used.

We now have that
\begin{equation}\label{hor}
{{\tau_Q}_*}_{|{T{\cal D}}}\circ {\bf h_{\Upsilon}} = 
{\bf h_{\gamma}}\circ {{\tau_Q}_*}_{|T{\cal D}}\;,
\end{equation}
i.e.\ the following diagram is commutative:
\[
\begin{array}{rcl}
(TTQ \supset)\; T{\cal D}& \stackrel{{\bf h_{\Upsilon}}}{\longrightarrow} 
& T{\cal D}\;(\subset TTQ)\\ 
{\tau_Q}_*\;\downarrow & & \downarrow\; {\tau_Q}_*\\
TQ & \stackrel{{\bf h_{\gamma}}}{\longrightarrow} & TQ
\end{array}
\]
Indeed, taking into account Proposition 3.1 we see that for any 
$w \in T_{v_q}{\cal D}$, $\gamma_q({\tau_Q}_*{\bf h_{\Upsilon}}(w)) = 
\Upsilon_{v_q}({\bf h_{\Upsilon}}(w)) = 0$ and, hence, 
${\tau_Q}_*{\bf h_{\Upsilon}}(w)$ is horizontal with respect to $\gamma$,
i.e.\ ${\tau_Q}_*{\bf h_{\Upsilon}}(w) \in {\cal D}_q$. 
By definition we also have ${\bf h_{\gamma}}({\tau_Q}_*w) \in {\cal D}_q$. 
Using the fact that $\pi \circ \tau_Q = \tau_{Q/G} \circ \pi_*$ we obtain:
\[
\pi_*({\bf h_{\gamma}}({\tau_Q}_*w)) = \pi_*{\tau_Q}_*w 
= {\tau_{Q/G}}_*(\pi_*)_*w = {\tau_{Q/G}}_*\rho_*w,
\]
where the last equation follows from (\ref{proj}).
Similarly, we have:
\[
\pi_*({\tau_Q}_*({\bf h_{\Upsilon}}w)) = {\tau_{Q/G}}_*(\pi_*)_*({\bf h_{\Upsilon}}w)
= {\tau_{Q/G}}_*\rho_*({\bf h_{\Upsilon}}w) = {\tau_{Q/G}}_*\rho_*w.
\]
We thus see that the $\gamma$-horizontal tangent vectors at $q$,
${\tau_Q}_*({\bf h_{\Upsilon}}w)$ and ${\bf h_{\gamma}}({\tau_Q}_*w)$, 
have the same projection under $\pi_*$ and, therefore, they are equal. 
This completes the proof of (\ref{hor}).
Denoting the curvature tensors of the principal connections $\gamma$ and
$\Upsilon$ by $\Omega^{\gamma}$ and $\Omega^{\Upsilon}$, respectively, one
can easily deduce from Proposition 3.1 and (\ref{hor}) the relation
\begin{equation}\label{curvature}
\Omega^{\Upsilon} = \tau_Q^*\Omega^{\gamma}.
\end{equation}
 
{\bf The 1-form.} Denote by $\theta'$ the pullback to ${\cal D}$ of the 
Poincar\'e-Cartan 1-form $\theta_L$, i.e.\ $\theta'=j_{\cal D}^* \theta_L$, 
where $j_{\cal D} : {\cal D} \hookrightarrow TQ$ is the canonical inclusion.
By means of the solution $X$ of (\ref{equat3}) we can construct a 1-form
$\alpha_X$ on ${\cal D}$ as follows:
\begin{equation}\label{alpha}
\alpha_X = i_X ({\bf h_{\Upsilon}}^*d\theta'-d{\bf h_{\Upsilon}}^*\theta') \; ,
\end{equation}
with the usual convention that, for an arbitrary $p$-form $\beta$ on ${\cal D}$, 
${\bf h_{\Upsilon}}^*\beta$ is the $p$-form defined by 
${\bf h_{\Upsilon}}^*\beta(X_1, \ldots, X_p) = 
\beta({\bf h_{\Upsilon}}(X_1), \ldots,{\bf h_{\Upsilon}}(X_p))$.

{\bf The Lagrangian.} The Lagrangian $L$ of the given mechanical system 
induces a Lagrangian $L^* : T(Q/G) \longrightarrow  \R$ on the quotient space 
$\bar{\cal D} \cong T(Q/G)$, given by $L^{*}(\bar{q},v_{\bar{q}}) = L(q,v^h_q)$
for any $q \in  \pi^{-1}(\bar{q})$, and where $v^h_q$ denotes the 
$\gamma$-horizontal lift of $v_{\bar{q}}$ at $q$. This is 
well-defined because of the $G$-invariance of $L$. Moreover, one can show 
that $L^*$ is a regular Lagrangian on $T(Q/G)$ (cf.\ \cite{LeMa}).

Now, we are in a position to state the following reduction result (see 
\cite{CaLeMaMa2,Ko}).

{\sc Proposition 3.2.} {\it The dynamics $X$ of the generalized Chaplygin system 
projects onto $\bar{{\cal D}}$, and its projection $\bar{X}$ is determined by
the equation 
\begin{equation}\label{redu}
i_{\bar{X}} \omega_{L^{*}} = dE_{L^{*}} + \overline{\alpha_X} \; ,
\end{equation}
where $\overline{\alpha_X}$ is the projection of the 1-form 
$\alpha_X$, defined by (\ref{alpha}). Moreover, we have that 
$i_{\bar{X}}\overline{\alpha_X}=0$.}

{\sc Remark 3.3.} The sign of $\overline{\alpha_X}$ in equation 
(\ref{redu}) differs from the one in the corresponding expression derived 
in \cite{CaLeMaMa2}, where the discussion took place in a more general
symplectic framework with an exact symplectic structure $\omega = d\theta$
(whereas here we have $\omega_L = - d\theta_L$). 
The signs would agree if we would have defined $\theta'$ as 
$- j_{\cal D}^* \theta_L$.

It can be easily verified that the form $\overline{\alpha_X}$ is a semi-basic
1-form on $\bar{\cal D}$ (see also Section 6), from which it then follows 
that the vector field $\bar{X}$, defined by (\ref{redu}), is a SODE. Moreover, one 
can show that not only the contraction of $\overline{\alpha_X}$ with 
$\bar{X}$ vanishes, but that, more generally, $i_Y\overline{\alpha_X} = 0$ 
for any SODE $Y$ on $\bar{\cal D}$. We thus see that a generalized Chaplygin 
system reduces to an unconstrained mechanical system, with an external 
nonconservative force of ``gyroscopic" type, which is geometrically represented by 
the 1-form $\overline{\alpha_X}$ (see also \cite{Ko,NeFu}). The 
``gyroscopic" character of this force is also in agreement with the fact 
that the projected energy function $E_{L^*}$ is a conserved quantity of the 
reduced dynamics. But there is more to be said about it.

{\sc Proposition 3.4.} {\it The 2-form 
$\Sigma={\bf h_{\Upsilon}}^*d\theta'-d{\bf h_{\Upsilon}}^*\theta'$ 
on ${\cal D}$ projects onto a 2-form $\bar{\Sigma}$ on $\bar{\cal D}$ and the 
1-form $\overline{\alpha_X}$ satisfies} 
\[
\overline{\alpha_X} = i_{\bar{X}} \bar{\Sigma} \, .
\]

{\it Proof.} Let $\xi_{\cal D}$ be the fundamental vector field of the 
$G$-action on ${\cal D}$, induced by an arbitrary element $\xi \in \hbox{\fr g}$. 
We must prove that $\xi_{\cal D}$ belongs to the characteristic distribution
of the 2-form $\Sigma$. First, we have that
\[
i_{\xi_{\cal D}} \Sigma = - i_{\xi_{\cal D}}d{\bf h_{\Upsilon}}^*\theta' \, .
\]
For any vector field $Y$ on ${\cal D}$, 
$i_{\xi_{\cal D}}d{\bf h_{\Upsilon}}^*\theta'(Y) = \xi_{\cal D}
(\theta' ({\bf h_{\Upsilon}}Y)) - \theta' ({\bf h_{\Upsilon}}[\xi_{\cal D},Y])$. 
Now, if $Y$ is vertical, we readily see that 
$i_{\xi_{\cal D}}d{\bf h_{\Upsilon}}^*\theta' (Y)=0$. If $Y$ is horizontal, we have
\[
i_{\xi_{\cal D}}d{\bf h_{\Upsilon}}^*\theta' (Y) = 
{\cal L}_{\xi_{\cal D}} \theta' (Y) = 0 \, ,
\]
because of the $G$-invariance of $\theta'$.
It therefore remains to prove that $i_{\xi_{\cal D}}d\Sigma =0$. 
For any two vector fields $Y,Z$ we have that
\begin{eqnarray*}
i_{\xi_{\cal D}}d\Sigma (Y,Z) &=& 
\left(i_{\xi_{\cal D}}d{\bf h_{\Upsilon}}^*d\theta'\right)(Y,Z)\\ 
&=& \xi_{\cal D}\left({\bf h_{\Upsilon}}^*d\theta'\right)(Y,Z) - 
{\bf h_{\Upsilon}}^*d\theta' ([\xi_{\cal D},Y],Z) + 
{\bf h_{\Upsilon}}^*d\theta' ([\xi_{\cal D},Z],Y) \, .
\end{eqnarray*}
If at least one of the vector fields $Y$ and $Z$ is vertical, then 
$i_{\xi_{\cal D}}d\Sigma (Y,Z)=0$. 
Taking $Y$ and $Z$ both horizontal, we find, taking into account
the $G$-invariance of $d\theta'$ and $\xi_{\cal D} \in {\frak X} ({\cal D})$,
\[
i_{\xi_{\cal D}}d\Sigma (Y,Z) = \xi_{\cal D}\left(d\theta'(Y,Z)\right) - 
d\theta' ([\xi_{\cal D},Y],Z) - d\theta' (Y,[\xi_{\cal D},Z]) = 
({\cal L}_{\xi_{\cal D}}d\theta')(Y,Z) = 0 \, .
\]
The last part of the proposition now immediately follows from (\ref{alpha})
and the projectability of $X$.
\QED

Consider a local trivialization $U \times G$ of $\pi$, with coordinates 
$(q^a,g^i)$, where $i = 1, \ldots k = \dim G$ and $a = 1, \ldots, n-k$.
Choosing a basis $e_i$ ($i = 1, \ldots, k$) of the Lie algebra $\hbox{\fr g}$,
and using the left trivialization $TG \cong G \times \hbox{\fr g}$, a tangent 
vector $v \in T_{(x,g)}(U \times G) \cong T_xU \times \hbox{\fr g}$ can 
be represented by a pair $(w,\xi)$, whereby $w \in T_xU$ and 
$\xi = \xi^ie_i \in \hbox{\fr g}$. In terms of the coordinates 
$(q^a, g^i, \dot{q}^a,\xi^i)$ on $T(U \times G)$ the $G$-invariant Lagrangian 
can then be written as 
\[
L= \ell(q^a,\dot{q}^a,\xi^i).
\]
Strictly speaking, $\ell$ represents the reduction of $L$ to $TQ/G$.
With respect to the given local trivialization, we further denote the connection 
coefficients of the given principal connection $\gamma$ by 
$\Gamma^i_a=\Gamma^i_a(q^1,\ldots,q^{n-k})$, and the constraints then take the form 
$\xi^i = -\Gamma^i_a\dot{q}^a$. In particular, it follows that the reduced 
Lagrangian $L^*$ is given by 
$L^*(q^a,\dot{q}^a) = \ell(q^a,\dot{q}^a,-\Gamma^i_b\dot{q}^b)$.

With all the above one can now derive the following coordinate expression for 
the reduced dynamics (see also \cite{BlKrMaMu,Ko}):
\begin{equation}\label{reduced}
\bar{X} = \dot{q}^a \frac{\partial}{\partial q^a} + \hat{W}^{ab} \left( \frac{\partial L^*}{\partial q^b} - \dot{q}^c \frac{\partial \hat{p}_b}{\partial q^c} - 
\alpha_b \right) \frac{\partial}{\partial \dot{q}^a} \, ,
\end{equation}
where $(\hat{W}^{ab})$ is the inverse of the Hessian matrix
$\displaystyle{\left(\frac{\partial^2 L^*}{\partial \dot{q}^a \partial \dot{q}^b}\right)}$, $\hat{p}_a = \displaystyle{\frac{\partial L^*}{\partial \dot{q}^a}}$,
and $\alpha_bdq^b$ is the local expression for the gyroscopic 1-form 
$\overline{\alpha_X}$. The $\alpha_b$ are explicitly given by 
\[
\alpha_b = -\left(\frac{\partial \ell}{\partial \xi^i}\right)^*
\left(\frac{\partial \Gamma^i_c}{\partial q^b}-\frac{\partial \Gamma^i_b}{\partial q^c}
- c^i_{jk}\Gamma^j_b\Gamma^k_c \right)\dot{q}^c \; ,
\]
where the * on the right-hand side indicates that, after computing the derivative
of $\ell$ with respect to $\xi^i$, one replaces  the $\xi^j$ everywhere by 
$-\Gamma^j_a\dot{q}^a$. The constants $c^i_{jk}$ appearing in the last term 
on the right-hand side are the structure constants of $\hbox{\fr g}$ with respect to 
the chosen basis, i.e.\ $[e_j,e_k] = c^i_{jk}e_i$. Note in passing that the
expressions $\displaystyle{\frac{\partial \Gamma^i_c}{\partial q^b}-
\frac{\partial \Gamma^i_b}{\partial q^c} - c^i_{jk}\Gamma^j_b\Gamma^k_c}$  
are the coefficients of the curvature of $\gamma$ in local form.

3.2. {\it Reconstruction}

A natural problem related to the reduction of mechanical systems with symmetry 
concerns the reverse procedure: once the solutions of the reduced dynamics 
have been obtained, how can one recover from it the solutions of the original 
system. This is called the ``reconstruction problem" of the dynamics. This 
problem is intimately related to the concepts of geometric and dynamic phase, 
which play an important role in various aspects of mechanics \cite{MaMoRa} and 
in the study of locomotion systems (for example, in the generation of net motion by 
cyclic changes in shape space \cite{KeMu,ShWi}). 

In the case of a generalized Chaplygin system, the reconstruction of the 
dynamics on ${\cal D}$ from the reduced dynamics on $\bar{\cal D}$ proceeds 
as follows \cite{CoLe}. Suppose that the flow of the reduced system $\bar{X}$ 
is known. Let $\bar{c}(t)$ be the integral curve of $\bar{X}$ starting 
at a given point $\bar{x} \in \bar{\cal D}$, and fix $x \in
\rho^{-1}(\bar{x})$. We want to find the integral curve $c(t)$ of $X$ with 
initial point $x$ and projecting onto $\bar{c}(t)$, i.e. $\rho(c(t)) = \bar{c}(t)$. 
But this is precisely the horizontal lift through $x$ of $\bar{c}(t)$, 
with respect to the principal connection $\Upsilon$. We recall here the proof 
of this simple fact (see also \cite{CoLe}). 

{\sc Proposition 3.5.} {\it The integral curve $c(t)$ of the generalized 
Chaplygin system, starting at $x \in {\cal D}$, is the horizontal lift with 
respect to the principal connection $\Upsilon$ of the integral curve $\bar{c}(t)$ 
of the reduced system starting at $\bar{x} = \rho(x)$.}

{\it Proof.} Let $d(t)$ denote the horizontal lift of $\bar{c}(t)$, starting at
$x$ (say, at $t=0$). Therefore, $\dot{d}(t) \in {\cal D}$ and 
$\rho(d(t)) = \bar{c}(t)$ for each $t$, and $d(0) = x$. 
Since $X$ and $\bar{X}$ are $\rho$-related, we have that $T\rho(X(d(t))) =
\bar{X}(\rho(d(t))) = \bar{X}(\bar{c}(t)) = T\rho(\dot{d}(t))$. 
Hence, $\dot{d}(t) - X(d(t))$ is vertical for each $t$. On the other hand, 
it is also horizontal, since $X \in {\cal U}$. Therefore, we deduce that 
$\dot{d}(t) = X(d(t))$. 
\QED

Summarizing: in the symplectic formalism, the reconstruction simply consists of
a horizontal lift operation with respect to the induced connection 
$\Upsilon$ living on ${\cal D}$.

\begin{center}
4. {\it Affine connection approach}
\end{center}

In this section we will describe the reduction and reconstruction problem
for generalized Chaplygin systems of mechanical type, from a different point
of view, namely in terms of the theory of affine connections.
First, we recall some general results which apply to any kind of
nonholonomic system of mechanical type. 

Let $Q$ be the configuration space of a mechanical system with Lagrangian 
of the form 
\[
L(v)=\frac{1}{2}{g}(v,v) - V \circ \tau_Q (v), \ v\in T_xQ \, ,
\]
where $g$ is a Riemannian metric on $Q$ and $V:Q \longrightarrow \R$ is the 
potential energy function. We denote by $\nabla^{g}$ the (covariant derivative 
operator of the) Levi-Civita connection associated to the metric 
${g}$. It is well-known that a curve $c: I \longrightarrow Q$ is a solution of 
the Euler-Lagrange equations for the (unconstrained) Lagrangian
if and only if
\[
\nabla^{g}_{\dot{c}(t)}\dot{c}(t) = - \hbox{grad}\,V(c(t)) \, ,
\]
where the gradient is also considered with respect to the metric $g$.

For nonholonomic systems there is a similar description. Let ${\cal D}$ again 
denote a (nonintegrable) distribution on $Q$, describing some linear nonholonomic 
constraints. The second-order differential equations (\ref{noholo}) for the 
mechanical system with Lagrangian $L$ and constraints ${\cal D}$, can be
written intrinsically as
\[
\nabla^{g}_{\dot{c}(t)}\dot{c}(t) + \hbox{grad}\,V(c(t)) \in 
{\cal D}_{\dot{c}(t)}^{\perp} \, , \; \; \dot{c}(t)\in {\cal D}_{c(t)} \, ,
\]
where ${\cal D}^{\perp}$ here denotes the ${g}$-orthogonal complement to ${\cal D}$
(see e.g.\ \cite{BlCr,Le,Sy,VeFa}).

Alternatively, if we denote by 
\[
{\cal P}: TQ\longrightarrow {\cal D}\, , \qquad 
{\cal Q}: TQ\longrightarrow {\cal D}^{\perp} \, ,
\]
the complementary ${g}$-orthogonal projectors, we can define an affine
connection 
\[
\bar{\nabla}_X Y=\nabla^{g}_X Y+(\nabla^{g}_X {\cal Q})(Y) \, ,
\]
such that the nonholonomic equations of motion can be rewritten as 
\[
\bar{\nabla}_{\dot{c}(t)}\dot{c}(t)= - {\cal P} (\hbox{grad}\,V(c(t))) \, ,
\]
and where we select the initial velocity  in ${\cal D}$ (cf. \cite{Le}). 

In what follows, we shall restrict our attention to Lagrangians of ``pure kinetic 
energy type", i.e. we assume $V = 0$. The reason for doing this is twofold. 
First, it makes the geometric picture more clear and tractable, in that the 
equations of motion for the nonholonomic mechanical system can then be seen 
as the geodesic equations of an affine connection. Secondly, the extention to  
systems with a nontrivial potential energy function is rather straightforward but, 
at least for those aspects of nonholonomic dynamics that are of interest to us 
here, it does not really tell us anything new.

It can be easily deduced from its definition that the connection $\bar{\nabla}$ 
restricts to ${\cal D}$, that is,  
\[
\bar{\nabla}_X Y\in {\cal D} \, ,
\]
for all $Y\in {\cal D}$ and $X\in {\frak X}(Q)$. This kind of affine
connections, which restrict to a given distribution, have been studied in
\cite{Le}. In particular, such a behaviour implies that the distribution
${\cal D}$ is {\bf geodesically invariant}, that is, for every geodesic
$c(t)$ of $\bar{\nabla}$ starting from a point in ${\cal D}$,
$\dot{c}(0) \in {\cal D}_{c(0)}$, we have that $\dot{c}(t) \in {\cal
D}_{c(t)}$. In \cite{Le}, a nice property is derived which
characterizes geodesically invariant distributions in terms of the
so-called symmetric product of vector fields, which is defined as
$\langle X:Y\rangle=\bar{\nabla}_XY + \bar{\nabla}_YX$. This property asserts that
${\cal D}$ is geodesically invariant if and only if we have that 
$\langle X:Y\rangle \in {\cal D}$, $\forall X$, $Y \in {\cal D}$. Note in passing 
that the symmetric product of vector fields is a differential geometric concept 
with important applications to control theory, first introduced 
in \cite{LewMur}. For instance, it plays a fundamental role in the controllability 
analysis of mechanical control systems \cite{Bullo2,CoMar,LewMur}, in the 
description of the evolution of these systems when starting from rest \cite{Bullo1}, 
and in the design of motion planning control algorithms \cite{Bullo2}.

{\sc Proposition 4.1.} {\it For all  $Z \in {\frak X}({Q})$ and $X$, 
$Y\in {\cal D}$ we have that}
\[
Z\left ({g}(X,Y)\right)={g}(\bar{\nabla}_Z X, Y)+{g}(X, \bar{\nabla}_Z Y) \, .
\]

{\it Proof.} In view of the definition of $\bar{\nabla}$, we have that
\begin{eqnarray*}
{g}(\bar{\nabla}_Z X, Y)+{g}(X, \bar{\nabla}_Z Y)&=& 
{g}(\nabla^{g}_Z X, Y)+{g}(X, \nabla^{g}_Z Y)+
{g}((\nabla^{g}_Z {\cal Q})(X), Y)+{g}(X, (\nabla^{g}_Z {\cal Q})(Y))\\
&=&Z\left ({g}(X,Y)\right) \, ,
\end{eqnarray*}
since $(\nabla^{g}_Z {\cal Q})(X),(\nabla^{g}_Z {\cal Q})(Y)\in {\cal
D}^{\perp}$ (see Proposition 6.1 in \cite{Le}). 
\QED

We derive from this proposition that the connection $\bar{\nabla}$ also has 
the following property: parallel transport is an isometry along the
distribution ${\cal D}$.  

A direct computation shows that the torsion of $\bar{\nabla}$ is the
skew-symmetric $(1,2)$-tensor field 
\[
\bar{T}(X,Y)=(\nabla^{g}_X {\cal Q})(Y)-(\nabla^{g}_Y {\cal Q})(X) \, .
\]
Observe that if $X, Y\in {\cal D}$ then $\bar{T}(X,Y)\in {\cal D}^{\perp}$. 

4.1. {\it Reduction}

We now return to the case of generalized Chaplygin systems of mechanical
type, but again assuming, for simplicity, that $L$ only consists of a 
kinetic energy part (i.e.\ $V=0$). From the above we then know that the equations
of motion of the system can be written as
\[
\bar{\nabla}_{\dot{c}(t)}\dot{c}(t) = 0, \quad \dot{c}(0) \in {\cal D}.
\]
The following is mainly inspired on the reduction 
process of generalized (or non-abelian) Chaplygin systems as described 
in \cite{Ko}. 

Let us define a metric $\tilde{{g}}$ on the base manifold $M (= Q/G)$ as 
follows
\[
\tilde{g}_x(u_x,v_x)={g}_q(U_q, V_q), \quad x \in M,\; u_x, v_x\in T_xM \, ,
\]
where $q\in \pi^{-1}(x)$ and $U_q, V_q$ are horizontal vectors which
project under $\pi$ onto $u_x$ and $v_x$, respectively. From the
$G$-invariance of $g$ we deduce that the right-hand side is independent of
the chosen point $q$ in the fibre $\pi^{-1}(x)$ and, hence, $\tilde{g}$ is
well defined. 

{\sc Proposition 4.2.} {\it We have that for all $X,Y\in {\frak X}(M)$ and 
$\xi\in \hbox{\fr g}$}
\[
{\cal L}_{\xi_Q}(\bar{\nabla}_{X^h}Y^h)=0 \, .
\]

{\it Proof.} Since $\xi_Q$ is a Killing vector field, it follows from Proposition
VI. 2.2 in \cite{KoNo} that 
\begin{equation}\label{asdf}
{\cal L}_{\xi_Q}\nabla^{g}_{X^h}Y^h=\nabla^{g}_{[\xi_Q,
X^h]}Y^h+\nabla^{g}_{X^h}[\xi_Q, Y^h]=0 \, , 
\end{equation}
because $Y^h$ and $X^h$ are projectable. Therefore, we only need to 
prove that
\[
{\cal L}_{\xi_Q}\left[\left(\nabla^{g}_{X^h}{\cal Q}\right)Y^h\right]=0 \, .
\]
This condition is equivalent to 
\begin{equation}\label{assi}
{\cal L}_{\xi_Q}\left[{\cal Q}(\nabla^{g}_{X^h}Y^h)\right]=0 \, .
\end{equation}
But, as ${\cal D}$ and ${\cal D}^{\perp}$ are $G$-invariant, we have 
that ${\cal L}_{\xi_Q}{\cal Q}=0$. This, together with (\ref{asdf}),
imply (\ref{assi}).
\QED

Now, we define an affine connection on $M$ as follows: 
for $X,Y \in {\frak X}(M)$, put
\[
\tilde{\nabla}_X Y= \pi_*(\bar{\nabla}_{X^h} Y^h)\, .
\]
This is well-defined since, by Proposition 4.2, the vector field 
$\bar{\nabla}_{X^h}Y^h$ is projectable, and one easily verifies that 
$\tilde{\nabla}$ satisfies the properties of an affine connection. 
Then, we obtain the following important result.

{\sc Proposition 4.3.}(\cite{Ko}) {\it The geodesics of $\bar{\nabla}$, 
with initial condition in ${\cal D}$, project onto the geodesics of $\tilde{\nabla}$.} 

{\it Proof.} Key fact for the proof is that ${\cal D}$ is geodesically invariant 
with respect to $\bar{\nabla}$.
\QED

Consequently, we have found that the equations of motion of the given 
generalized Chaplygin system reduce to the geodesic equations of the
induced affine connection $\tilde{\nabla}$ on $M = Q/G$. 

Consider the following (0,3)-tensor field on $Q$:
\[
K_q (U_q,V_q,W_q) = g_q({\bf h_{\gamma}}U_q,(\Omega^{\gamma}(V_q,W_q))_Q(q)) \, ,
\]
where ${\bf h_{\gamma}}$ is the horizontal projector and $\Omega^{\gamma}$ is 
the curvature of the connection $\gamma$. Observe that $K$ is horizontal, 
i.e.\ it vanishes if one of its arguments is a vertical vector, and it is 
skew-symmetric in its last two arguments. Moreover, one can see that
\[
K_{gq} (T\Psi_g (U_q),T\Psi_g (V_q),T\Psi_g(W_q)) = K_q (U_q,V_q,W_q) \, ,
\]
for all $g \in G$ and $q \in Q$. Consequently, $K$ induces a (0,3)-tensor 
on the base manifold $M$
\[
\tilde{K}_x(u_x,v_x,w_x) = K_q (U_q,V_q,W_q) \, ,
\]
where $\pi(q)=x$ and $U_q,V_q,W_q$ are tangent vectors in $q$ projecting
onto $u_x,v_x,w_x$, respectively. $K$ (resp.\ $\tilde{K}$) is called the 
{\bf metric connection tensor} on $Q$ (resp.\ $M$).

In \cite{Ko}, it was shown that application of the so-called Hamel's approach
to mechanics, leads to two additional affine connections on $M$, whose 
geodesics are also solutions of the reduced nonholonomic problem. 
These connections are given by
\begin{equation}\label{hamel}
({\nabla^{H}_1})_X Y = \nabla^{\tilde{g}}_X Y + B(X,Y) \; , \quad  
({\nabla^{H}_2})_X Y = \nabla^{\tilde{g}}_X Y + B(Y,X) \, ,
\end{equation}
where $B$ is the (1,2)-tensor field defined by 
$\beta (B(X,Y)) = \tilde{K}(X,Y, \sharp_{\tilde{g}}\beta)$, for any 
$\beta \in \Lambda^1(Q/G)$, $X$, $Y \in {\frak X}(Q/G)$. In general, 
the tensor which encodes the difference between an affine connection $\nabla$ 
on a Riemannian manifold and the Levi-Civita connection corresponding to the 
Riemannian metric, is called the {\bf contorsion} of $\nabla$. So, 
$B$ here represents the contorsion of $\nabla^H_1$.

The following explicit formula for the connection $\tilde{\nabla}$ 
was then derived in \cite{Ko} (up to a minor misprint): 
\[
\tilde{\nabla}_X Y=\nabla^{\tilde{g}}_X Y + 
\frac{1}{2} \left( B(X,Y)+B(Y,X)-C(X,Y) \right) \, ,
\]
where $C$ is the (1,2)-tensor field implicitly defined by $\beta (C(X,Y)) = K(\sharp_{\tilde{g}}\beta,X,Y)$, for arbitrary $\beta \in \Lambda^1(Q/G)$, 
$X$, $Y \in {\frak X}(Q/G)$. As noted in \cite{Ko}, the average of 
Hamel's connections, i.e.\  
$\nabla^{H/2} = \frac{1}{2}(\nabla_{1}^H + \nabla_{2}^H)$, in general differs 
from $\tilde{\nabla}$, because of the skew-symmetric term $C(X,Y)$. 

It is interesting to observe that from Proposition 4.1 one can easily deduce
\[
\tilde{\nabla}\tilde{g}=0 \, ,
\]
that is, $\tilde{\nabla}$ is a {\bf metric connection}. From 
Proposition A.3 (see Appendix) and the definition of $B$, it is 
readily seen that $\nabla^H_1$ is also a metric connection. In general, however, 
the connections $\nabla^H_2$ and $\nabla^{H/2}$ will not be metric. In fact, 
it is not hard to prove the following result.

{\sc Proposition 4.4.} {\it The following properties are equivalent:
\begin{enumerate}
\item $\nabla^H_2$ is metric;
\item $\nabla^{H/2}$ is metric;
\item the tensor field $B$ is skew-symmetric;
\item $\nabla^{H/2}$ is the Levi-Civita connection of $\tilde{g}$.
\end{enumerate}}

Later we will see that these properties are also equivalent to the 
vanishing of the 1-form $\overline{\alpha_X}$ and, hence, to the Hamiltonian
nature of the reduced system (cf.\ Corollary 6.2).
 
The torsion of $\tilde{\nabla}$ is given by $\tilde{T}(X,Y) = 
\pi_*\,\bar{T}(X^h,Y^h)$, with $\bar{T}$ the torsion of $\bar{\nabla}$. 
Then, we see from the above that the metric connection $\tilde{\nabla}$ is the 
Levi-Civita connection associated to $\tilde{g}$ iff the torsion of $\bar{\nabla}$ 
takes values in the vertical tangent bundle to $\pi$ for each pair of vectors in 
${\cal D}$.

Finally, the following result shows that equality of $\tilde{\nabla}$ and 
$\nabla^{H/2}$ is a rather strong condition.

{\sc Proposition 4.5.}
\[
\tilde{\nabla} = \nabla^{H/2} \; \Longleftrightarrow \; 
\tilde{\nabla} = \nabla^{H/2} = \nabla^{\tilde{g}} \, .
\]

{\it Proof.} If $\tilde{\nabla} = \nabla^{H/2}$, then $\nabla^{H/2}$ is metric. 
By Proposition 4.4, this implies that $\nabla^{H/2}$ coincides 
with $\nabla^{\tilde{g}}$. The reverse implication is trivial.
\QED

4.2. {\it Reconstruction} 

In view of Proposition 4.3 above, we see that, in the present
setting, the reconstruction of the solution curves in $Q$ of the given
constrained system, from those of the reduced system on $M$, consists of 
a horizontal lift operation with respect to the connection $\gamma$. That is, 
let $\tilde{c}(t)$ be a geodesic of $\tilde{\nabla}$ and choose 
$c(0) \in Q$ such that $\pi(c(0))=\tilde{c}(0)$. 
Then, the geodesic starting at $c(0)$, with initial velocity 
$\dot{\tilde{c}}(0) \in {\cal D}_{c(0)}$, is precisely the horizontal lift 
of $\tilde{c}(t)$ with respect to the principal connection $\gamma$.

\begin{center}
5. {\it Examples}
\end{center}

5.1. {\it Mobile robot with fixed orientation} 

Consider the motion of a mobile robot whose body maintains a fixed 
orientation with respect to its environment (see \cite{KeMu} for more details).
The robot has three wheels, with radius $R$, which turn simultaneously
about independent axes, and perform a rolling without sliding over a 
horizontal floor. Let $(x,y) \in \R^2$ denote the position of the center 
of mass of the robot (in a cartesian reference frame, with horizontal $x$- and 
$y$-axis), $\theta \in \S^1$ the steering angle of the wheels and $\psi \in \S^1$ 
the rotation angle of the wheels in their rolling motion over the floor. 
The configuration space can then be modelled by $Q=\S^1 \times \S^1 \times \R^2$.

The Lagrangian function $L$ is the kinetic energy function corresponding to 
the metric
\[
g=mdx \otimes dx + m dy \otimes dy + J d\theta \otimes d\theta + 
3 J_w d\psi \otimes d\psi \, ,
\]
where $m$ is the mass of the robot, $J$ its moment of inertia and $J_w$ the 
axial moment of inertia of each wheel. The constraints are induced by the conditions
that the wheels roll without sliding, in the direction in which they ``point",
and that the instantaneous contact points of the wheels with the floor have 
no velocity component orthogonal to that direction (cf.\ \cite{KeMu}):
\begin{eqnarray*}
\dot{x} \sin \theta - \dot{y} \cos \theta &=& 0 \, , \\
\dot{x} \cos \theta + \dot{y} \sin \theta - R \dot{\psi} &=& 0 \, .
\end{eqnarray*}
The constraint distribution ${\cal D}$ is then spanned by 
\[
\left\{ \frac{\partial}{\partial \theta} \, , \; 
\frac{\partial}{\partial \psi} + R\left( \cos \theta \frac{\partial}{\partial x} 
+ \sin \theta \frac{\partial}{\partial y} \right) \right\} \, .
\]

If we consider the abelian action of $G=\R^2$ on $Q$ by translations
\[
\begin{array}{rccl}
\Psi: & G \times Q & \longrightarrow & Q \\
& ((a,b),(\theta,\psi, x,y)) & \longmapsto & (\theta,\psi,a+x,b+y) \, ,
\end{array}
\]
we see that the constraint distribution ${\cal D}$ can be interpreted as 
the horizontal subspace of the principal connection 
$\gamma = (dx-R \cos \theta d\psi)e_1 + (dy-R \sin \theta d\psi)e_2$, 
where $\{ e_1, e_2 \}$ is the canonical basis of $\R^2$ (identified with
the Lie algebra of $G$).

The metric induced on $M=Q/G$ here becomes
\[
\tilde{g} = J d\theta \otimes d\theta + (3J_\omega + mR^2) d\psi \otimes d\psi \, .
\]
The reduced Lagrangian $L^*$ is the kinetic energy function corresponding to the 
metric $\tilde{g}$. Moreover, one easily verifies that the gyroscopic 1-form 
$\overline{\alpha_X}$ identically vanishes and, hence, the symplectic 
reduction (\ref{redu}) yields
\[
i_{\bar{X}} \omega_{L^{*}} = dE_{L^{*}} \, ,
\]
i.e.\ the reduced system is an unconstrained, purely Lagrangian system.

On the other hand, one can easily check that, in this example, the metric 
connection tensor $\tilde{K}$ also vanishes. Consequently, with the notations of 
the previous section, $\tilde{\nabla}=\nabla^{H/2}$ and by Proposition 4.5 
we then have that $\tilde{\nabla}=\nabla^{H/2} = \nabla^{\tilde{g}}$.

5.2. {\it Two-wheeled planar mobile robot}

Consider the motion of a two-wheeled planar mobile robot (or ``two-wheeled 
carriage") which is able to move in the direction in which it points 
and, in addition, can spin about a vertical axis \cite{KeMu,Ko,LeMa,NeFu}. 
Let $P$ be the intersection point of the horizontal symmetry axis of the
robot and the horizontal line connecting the centers of the 
two wheels. The position and orientation of the robot is determined, with respect to a 
fixed cartesian reference frame (with horizontal $x$- and $y$-axes) 
by $(x,y,\theta) \in SE(2)$, where $\theta \in \S^1$ is the heading angle 
and the coordinates $(x,y) \in \R^2$ locate the point $P$ (see \cite{Ko,NeFu}). 
Let $\psi_1$, $\psi_2 \in \S^1$ denote the rotation 
angles of the wheels which are assumed to be controlled independently and roll
without slipping on the floor. The configuration space of this system is 
$Q=\S^1 \times \S^1 \times SE(2)$.

The Lagrangian function is the kinetic energy corresponding to the metric
\begin{eqnarray*}
g&=&m dx \otimes dx + m dy\otimes dy +  
m_0 l \cos\theta (dy \otimes d\theta + d\theta \otimes dy) \\
&&- m_0l\sin\theta (dx\otimes d\theta + d\theta \otimes dx)+ Jd\theta\otimes d\theta + J_2d\psi_1\otimes d\psi_1+J_2 d\psi_2\otimes d\psi_2 \, ,
\end{eqnarray*}
where $m=m_0 + 2m_1$, $m_0$ is the mass of the robot without the wheels, 
$J$ its moment of inertia with respect to the vertical axis, 
$m_1$ the mass of each wheel, $J_w$ the axial moments of inertia of the wheels, 
and $l$ the distance between the center of mass $C$ of the robot and the point $P$.

The constraints, induced by the conditions that there is no lateral sliding of
the robot and that the motion of the wheels also consists of a rolling without 
sliding are:
\begin{eqnarray*}
\dot x\sin\theta-\dot y \cos\theta&=&0 \, , \\
\dot x\cos\theta+\dot y\sin\theta+c\dot\theta+R \dot{\psi}_1&=&0 \, , \\
\dot x\cos\theta+\dot y\sin\theta-c\dot\theta+R \dot{\psi}_2&=&0\, ,
\end{eqnarray*}
where $R$ is the radius of the wheels and $2c$ the lateral length of the robot. 
The constraint distribution ${\cal D}$ is then spanned by 
\[
\left\{ \frac{\partial}{\partial \psi_1} - \frac{R}{2} \left( \cos \theta \frac{\partial}{\partial x} + \sin \theta \frac{\partial}{\partial y} + \frac{1}{c} \frac{\partial}{\partial \theta} \right) \, , \; 
\frac{\partial}{\partial \psi_2} - \frac{R}{2} \left( \cos \theta \frac{\partial}{\partial x} + \sin \theta \frac{\partial}{\partial y} - \frac{1}{c} \frac{\partial}{\partial \theta} \right)
\right\} \, .
\]

If we consider the action of $G=SE(2)$ on $Q$
\[
\begin{array}{rccl}
\Psi: & G \times Q & \longrightarrow & Q \\
&((a,b,\alpha),(\psi_1,\psi_2,x,y,\theta)) & \longmapsto & (\psi_1,\psi_2,
a+x\cos\alpha-y\sin\alpha, b+x\sin\alpha+y\cos\alpha,\alpha+\theta ) \, ,
\end{array}
\]
we see that the constraint distribution ${\cal D}$ can be interpreted as the 
horizontal subspace of the principal connection
\begin{eqnarray*}
\gamma &=& \left(dx + \frac{R}{2} \cos \theta d\psi_1 + \frac{R}{2} \cos \theta d\psi_2 + y(d\theta + \frac{R}{2c} d\psi_1 - \frac{R}{2c} d\psi_2)\right)e_1 \\
&& + \left(dy + \frac{R}{2} \sin \theta d\psi_1 
+ \frac{R}{2} \sin \theta d\psi_2-x(d\theta + \frac{R}{2c} d\psi_1 - \frac{R}{2c} d\psi_2)\right) e_2 \\
&& + (d\theta + \frac{R}{2c} d\psi_1 - \frac{R}{2c} d\psi_2)e_3 \, ,
\end{eqnarray*}
where $\{ e_1, e_2, e_3 \}$ is the canonical basis of the Lie algebra of $G$, with 
associated fundamental vector fields
\[
(e_1)_Q = \frac{\partial }{\partial x} \, , \; \; 
(e_2)_Q = \frac{\partial }{\partial y} \, , \; \;
(e_3)_Q = \frac{\partial }{\partial \theta} -y \frac{\partial }{\partial x} + 
 x \frac{\partial }{\partial y} \, .
\]
The curvature of $\gamma$ is
\[
\Omega= \frac{R^2}{2c}(\sin \theta \, e_1 - \cos\theta \, e_2) d\psi_1 \wedge d\psi_2  \, .
\]

The induced metric on $M=Q/G$ is given by
\begin{eqnarray*}
\tilde{g} &=& (J_2 + m \frac{R^2}{4} + J \frac{R^2}{4c^2}) d\psi_1 \otimes d\psi_1 \\
&& + (m \frac{R^2}{4} - J \frac{R^2}{4c^2}) (d\psi_1 \otimes d\psi_2 
+ d\psi_2 \otimes d\psi_1) + (J_2 + m \frac{R^2}{4} + J \frac{R^2}{4c^2}) 
d\psi_2 \otimes d\psi_2 \, .
\end{eqnarray*}
The Lagrangian $L^*$ is the kinetic energy function induced by $\tilde{g}$. The gyroscopic 
1-form $\overline{\alpha_X}$ here becomes
\[
\overline{\alpha_X} = \frac{m_0 l R^3}{4c^2} (\dot{\psi_2}-\dot{\psi_1}) 
(\dot{\psi_1}d\psi_2-\dot{\psi_2}d\psi_1) \, .
\]
Then, the symplectic reduction (\ref{redu}) yields 
\[
i_{\bar{X}} \omega_{L^{*}} = dE_{L^{*}} + \overline{\alpha_X} \, .
\]
On the other hand, the metric connection tensor $\tilde{K}$ is given by
\[
\tilde{K}= \frac{m_0 l R^3}{4c^2} (d\psi_1 \otimes d\psi_1 \wedge d\psi_2 - d\psi_2 \otimes d\psi_1 \wedge d\psi_2) \, .
\]
It is easily seen that, in this case, the tensor field $B$ is not skewsymmetric and so, 
by Proposition 4.4, $\nabla^{H/2} \not = \nabla^{\tilde{g}}$. In addition, the Christoffel symbols of the metric connection $\tilde{\nabla}$ are
\[
\begin{array}{rcl}
\Gamma_{\psi_1 \psi_1}^{\psi_1} &=& K_1 \, , \\
\Gamma_{\psi_1 \psi_2}^{\psi_1} &=& -K_2 \, ,  \\
\Gamma_{\psi_2 \psi_1}^{\psi_1} &=& -K_1 \, ,  \\
\Gamma_{\psi_2 \psi_2}^{\psi_1} &=& K_2 \, , \, ,
\end{array}
\begin{array}{rcl}
\Gamma_{\psi_1 \psi_1}^{\psi_2} &=& K_2 \, , \\
\Gamma_{\psi_1 \psi_2}^{\psi_2} &=& -K_1 \, , \\
\Gamma_{\psi_2 \psi_1}^{\psi_2} &=& -K_2 \, , \\
\Gamma_{\psi_2 \psi_2}^{\psi_2} &=& K_1 \, ,
\end{array}
\]
where
\begin{eqnarray*}
K_1 &=& m_0 R^5 l \frac{J_2-mc^2}{4c^2(2Jc^2+j_2R^2)(2J+mc^2)} \, , \\
K_2 &=& R^3 l \frac{4Jc^2+(J_2+mc^2)R^2m_0}{4c^2(2Jc^2+j_2R^2)(2J+mc^2)} \, .
\end{eqnarray*}
Clearly, the torsion $\tilde{T}$ does not vanish, and so $\tilde{\nabla} 
\not =\nabla^{\tilde{g}}$. 
\newpage
\begin{center}
6. {\it Relation between both approaches}
\end{center}

The above examples show us the following intriguing fact. In the case
of the mobile robot with fixed orientation, the 1-form $\overline{\alpha_X}$
identically vanishes, and thus the reduced problem has no external gyroscopic
force in its unconstrained symplectic formulation. Consequently, since there
is no potential, the solutions of the reduced system are geodesics of the 
Levi-Civita connection $\nabla^{\tilde{g}}$. Indeed, we verified that $\tilde{\nabla}=\nabla^{H/2} = \nabla^{\tilde{g}}$. 
However, in the case of the two-wheeled mobile robot we obtained 
$\overline{\alpha_X} \not = 0$ and $\tilde{\nabla} \not = \nabla^{\tilde{g}} 
\not = \nabla^{H/2}$. Apparently, there exists a relation between the properties 
of the contorsions of the connections considered in Section 4 and the vanishing 
(or not) of the gyroscopic 1-form.

Using the definition of $\alpha_X$ in (\ref{alpha}), one can check that the
following relation holds (see \cite{CaLeMaMa2}): for any $Y \in {\frak X}({\cal D})$,
\[ 
\alpha_X(Y)= {\bf v_{\Upsilon}}(Y)(\theta_L(X))+\theta_L(R(X,Y))+
\theta_L({\bf h_{\Upsilon}}[X,{\bf v_{\Upsilon}}(Y)]) \, ,
\]
where $R$ is the tensor field of type (1,2) on ${\cal D}$ given by 
\[
R = \frac{1}{2} [{\bf h_{\Upsilon}}, {\bf h_{\Upsilon}}] \; ,
\]
with $[\;,\;]$ denoting the Nijenhuis bracket of type (1,1) tensor fields.
The relation between $R$ and $\Omega^{\Upsilon}$, the curvature tensor of the 
principal connection $\Upsilon$, is given by 
$R(U,V)(v) = ({\Omega}^{\Upsilon}(U_v,V_v))_{TQ}(v)$ for any 
$U,V \in {\frak X}({\cal D})$ and $v \in {\cal D}$.

In particular, if we take a horizontal vector field $Y \in {\cal U}$, we deduce 
from the above that $\alpha_X(Y)=\theta_L (R(X,Y))$.
Herewith, the action of the gyroscopic 1-form $\overline{\alpha_X}$ on a vector field 
$Z \in {\frak X}(T(Q/G))$, evaluated at a point 
$w_{\bar{q}} \in T(Q/G) (\cong \bar{\cal D})$, becomes  
\begin{eqnarray*}
\overline{\alpha_X}(Z)(w_{\bar{q}})
&=& {\alpha_X}(Z^h)(v_q) = ({\theta_L})_{v_q}(R(X,Z^h)(v_q)) \\
&=& ({\theta_L})_{v_q}\left((\Omega^{\Upsilon}(X_{v_q},Z^h_{v_q}))_{TQ}(v_q)\right)
=({\theta_L})_{v_q}\left((\Omega^{\gamma}({\tau_Q}_*X_{v_q},
{\tau_Q}_*Z^h_{v_q}))_{TQ}(v_q)\right)\;, 
\end{eqnarray*}
for an arbitrary $v_q \in {\cal D}$ such that $\rho(v_q)=w_{\bar{q}}$, and 
where the last equality has been derived using (\ref{curvature}). In these 
expressions, $Z^h$ is the horizontal lift of $Z$ with respect to $\Upsilon$. 
Recalling that the Poincar\'e-Cartan 1-form $\theta_L$ and the Legendre mapping 
(induced by the given Lagrangian) $FL : TQ \longrightarrow T^*Q$ are related by 
$({\theta_L})_{v_q}(u) = <FL(v_q),{\tau_Q}_*(u)>$, for any $u \in T_{v_q}TQ$,
and taking into account that $X$ is a SODE, we further obtain
\begin{eqnarray}\label{baralpha}
\overline{\alpha_X}(Z)(w_{\bar{q}})
&=& <FL(v_q), (\Omega^{\gamma}({\tau_Q}_*X_{v_q},{\tau_Q}_*Z^h_{v_q}))_Q(q)> \nonumber \\
&=& g_{q}\left(v_q,(\Omega^{\gamma}(v_q,{\tau_Q}_*Z^h_{v_q}))_Q(q)\right) \nonumber \\
&=& g_{q}\left(v_q,(\Omega^{\gamma}(v_q,({\tau_{Q/G}}_*Z_{w_{\bar{q}}})^h_q))_Q(q)\right)\;. 
\end{eqnarray}
Note that in the last expression, the horizontal lift of 
${\tau_{Q/G}}_*Z_{w_{\bar{q}}}$ is the one with respect to $\gamma$. 

An important observation is that (\ref{baralpha}) immediately shows that
the gyroscopic 1-form $\overline{\alpha_X}$ is semi-basic with respect to 
the canonical fibration $\tau_{Q/G}: T(Q/G) \longrightarrow Q/G$. Indeed,
assume $({\tau}_{Q/G})_*\circ Z = 0$, then it readily follows that
$\overline{\alpha_X}(Z) = 0$.

Elaborating (\ref{baralpha}) a bit further, using the metric connection tensor
$\tilde{K}$ and the contorsion $B$ introduced in Section 4.1, we find
\begin{eqnarray*}
\overline{\alpha_X}(Z)(w_{\bar{q}})
&=& g_q\left(v_q, \Omega(v_q,({\tau_{Q/G}}_*Z_{w_{\bar{q}}})^h_q))_Q(q)\right) \\ 
&=& \tilde{K}_{\bar{q}} (w_{\bar{q}},w_{\bar{q}},{\tau_{Q/G}}_*Z_{w_{\bar{q}}}) \\
&=& \tilde{g}_{\bar{q}} (B(w_{\bar{q}},w_{\bar{q}}),{\tau_{Q/G}}_*Z_{w_{\bar{q}}}) \, .
\end{eqnarray*}
This proves the next result, which was already implicit in the work of 
Koiller \cite{Ko}.

{\sc Proposition 6.1.} {\it An explicit relation between the gyroscopic 1-form 
and the contorsion tensor field $B$, defined in Section 4, is given by
\[
(\overline{\alpha_X})_{w_{\bar{q}}}(u)
= \tilde{g}_{\bar{q}}(B(w_{\bar{q}},w_{\bar{q}}),{\tau_{Q/G}}_*u) \, ,
\]
for all $u \in T_{w_{\bar{q}}}T(Q/G))$, $w_{\bar{q}} \in T(Q/G)$.}

From this we can immediately deduce, taking into account Proposition 4.4: 

{\sc Corollary 6.2.} {\it The following statements are equivalent:
\begin{enumerate}
\item $\overline{\alpha_X}$ vanishes identically;
\item $B$ is skew-symmetric;
\item $\nabla^{H/2} = \nabla^{\tilde{g}}$.
\end{enumerate}}

{\sc Remark 6.3.} One can think of simple examples in which $B$ is skew-symmetric 
but nonzero. Consequently, if $\overline{\alpha_X}$ vanishes, this does not imply 
$\tilde{\nabla} = \nabla^{\tilde{g}}$, although in such a case both connections 
do have the same geodesics (see Appendix).

By means of Proposition 6.1, one can also recover the gyroscopic 
character of $\overline{\alpha_X}$, established already in 
Proposition 3.4. For that purpose, let us define the following 2-form
on $T(Q/G)$:
\begin{equation}\label{dosforma}
\Xi(Y,Z)(w_{\bar{q}}) = \tilde{g}_{\bar{q}}(B(w_{\bar{q}},{\tau_{Q/G}}_*Y_{w_{\bar{q}}}),
{\tau_{Q/G}}_*Z_{w_{\bar{q}}}) \, .
\end{equation}
One readily verifies that $\Xi$ is indeed bilinear and by Proposition A.3,
$\Xi(Y,Y)=0$. It is then easy to check that 
\[
\overline{\alpha_X} = i_{\bar{X}} \Xi \, .
\]
In local coordinates $q^a\;(a= 1, \ldots, n-k)$ on $M = Q/G$, we have that
\begin{equation}\label{local}
\Xi = \sum_{a<b}\dot{q}^e B^c_{ea}\tilde{g}_{bc}dq^a \wedge  dq^b \, , \qquad
\overline{\alpha_X} = \sum_{a,e,c} \dot{q}^a \dot{q}^e B^c_{ea}\tilde{g}_{bc} dq^b \, .
\end{equation}

A careful calculation, very similar to the one performed for proving 
Proposition 6.1, reveals that for generalized Chaplygin systems
of mechanical type, the 2-form $\bar{\Sigma}$ of Proposition 3.4 
and the above 2-form $\Xi$ coincide.

\begin{center}
7. {\it Integrability aspects and the existence of an invariant measure}
\end{center}

7.1. {\it Koiller's question}

In \cite{Ko}, the author wonders whether there might exist an invariant measure 
for the reduced equations of a generalized Chaplygin system. In the following, we 
shall deal with this problem in some detail.

The existence of a measure which is invariant under the flow of a given 
dynamical system is a strong property. Indeed, using an integrating factor 
it is possible to derive from it (locally) an integral of the motion.
This fact plays an important role in discussions concerning the integrability 
of the system under consideration, as illustrated by the following theorem.

{\sc Theorem 7.1.}(\cite{Ar})
{\it Suppose that the system $\dot{x}=X(x)$, $x \in N$, with $N$ an $n$-dimensional
smooth manifold, admits an invariant measure and has $n-2$ first 
integrals $F_1,...,F_{n-2}$. Suppose also that $F_1,...,F_{n-2}$ are 
independent on the invariant set $N_c=\{ x\in N : F_s(x)=c_s, 1 \le s \le n-2 \}$. 
Then:
\begin{itemize}
\item the solutions of the differential equation lying on $N_c$ can be found 
by quadratures. 
\end{itemize}
Moreover, if $L_c$ is a compact connected component of the level 
set $N_c$ and if $X$ does not vanish on $L_c$, then
\begin{itemize}
\item $L_c$ is a smooth manifold diffeomorphic to a two-torus;
\item one can find angular coordinates $\varphi_1$, $\varphi_2 \,\hbox{mod} \,(2\pi)$ 
on $L_c$ in terms of which the differential equations take the simple form
\[ 
\dot{\varphi}_1 = \frac{\omega_1}{\Phi(\varphi_1,\varphi_2)} \, , \; \; 
\dot{\varphi}_2 = \frac{\omega_2}{\Phi(\varphi_1,\varphi_2)} \, ,
\]
where $\omega_1$, $\omega_2$ are constant and $\Phi$ is a smooth positive 
function which is $2\pi$-periodic in $\varphi_1$, $\varphi_2$.
\end{itemize}}

By the Riesz representation theorem, we know that each volume form on an 
orientable manifold induces a unique measure on the Borel 
$\sigma$-algebra \cite{AbMa}. Therefore, with a view on tackling the integrability
problem of generalized Chaplygin systems, it is worth looking for invariant 
volume forms under the flow of the reduced dynamics $\bar{X}$. This is what we
intend to do in the sequel.

In this section, we again consider Chaplygin systems of mechanical type
with a Lagrangian $L=\frac{1}{2}g - V$, where the metric $g$ and the 
potential energy function $V$ are $G$-invariant. The reduced 
equations of motion are (cf.\ (\ref{redu}))
\begin{equation}\label{reduKo}
i_{\bar{X}} \omega_{L^{*}} = dE_{L^{*}} + \overline{\alpha_X} \; ,
\end{equation}
where $L^*=\frac{1}{2}\tilde{g}-\tilde{V}$, with $\tilde{g}$ and $\tilde{V}$ 
the metric and the potential function on $M$ induced by, respectively, $g$ and 
$V$. Remember that the energy $E_{L^*}$ is a constant of the 
motion. The local expression for the reduced dynamics takes the form
(cf.\ (\ref{reduced})) 
\begin{eqnarray*}
\bar{X} = \dot{q}^a \frac{\partial}{\partial q^a} - \left( \tilde{g}^{ab} (\alpha_b + \frac{\partial \tilde{V}}{\partial q^b}) + \dot{q}^b \dot{q}^c \tilde{\Gamma}_{bc}^a \right) \frac{\partial}{\partial \dot{q}^a} \, ,
\end{eqnarray*}
where $\tilde{\Gamma}_{bc}^a$ are the Christoffel symbols of the Levi-Civita 
connection $\nabla^{\tilde{g}}$.

The gyroscopic systems usually encountered in the mechanics literature 
\cite{WaKr,Ya} differ in a crucial way from the ones we obtain through the reduction 
of a generalized Chaplygin system. In fact, the common situation in mechanics is that 
of a system, with configuration space $P$, described by an equation of the form  
\[
i_{\Gamma} \omega_{\L} = dE_{\L} + \alpha \, ,
\]
where ${\L}:TP \longrightarrow \R$ is a (regular) Lagrangian and where the gyroscopic
force is represented by a 1-form $\alpha = i_{\Gamma}(\tau_P^* \Pi)$, with $\Pi$ a 
closed 2-form on $P$. These systems are then Hamiltonian with respect to the 
symplectic 2-form $\omega = \omega_{\L} - \tau_P^*\Pi$, and thus they admit an 
invariant measure, determined by the volume form $\omega^n= \omega_{\L}^n$. 

In some sense, our reduced system (\ref{reduKo}) exhibits the opposite behaviour. 
Indeed, the 2-form $\Xi$, defined by (\ref{dosforma}), is semi-basic but in general
it is not basic, i.e.\ it is not the pull-back of a 2-form on the base $M = Q/G$. 
This can be readily deduced from its local expression (cf.\ \ref{local}). 
Moreover, using (\ref{local}), the following property is easily proved.

{\sc Proposition 7.2.} {\it The 2-form $\Xi$ is closed if and only if it is 
identically zero.}

Note, in passing, that a similar property also applies to the gyroscopic 1-form
$\overline{\alpha_X}$. The semi-basic character of $\Xi$ ensures, however, that 
the 2-form $\omega_{L^*} - \Xi$ is still nondegenerate and, consequently, we have 
that the equation (\ref{reduKo}) can be rewritten in the form
\begin{equation}\label{as}
i_{\bar{X}} \omega = dE_{L^*} \, ,
\end{equation}
with $\omega=\omega_{L^*} - \Xi$ an almost symplectic form (i.e.\ a nondegenerate,
but not necessarily closed 2-form).

In \cite{St}, S.V. Stanchenko has studied Chaplygin systems of mechanical type 
with an abelian Lie group in terms of differential forms, in a way which shows 
many links to the approach described in Section 3. In our setting, his 
results can be generalized to the non-abelian case for any kind of 
generalized Chaplygin system.

Let us assume, following \cite{St}, that there exists a function 
$F \in C^{\infty}(T(Q/G))$ such that
\begin{equation}\label{as2}
dF \wedge \theta_{L^*}= \Xi \, .
\end{equation}
Putting $N=\exp F$, we have that
\[
d(N\omega)=d(N\omega_{L^*}-N\Xi)=d(N\omega_{L^*}-dN \wedge \theta_{L^*}) = 0 \, .
\]
Since (\ref{as}) can still be written as
\[
i_{\bar{X}/N} (N\omega) = dE_{L^*} \, ,
\]
we deduce that ${\cal L}_{\bar{X}/N}(N\omega) = 0$. Consequently,
\[
0 = {\cal L}_{\bar{X}/N}(N\omega)^n = {\cal L}_{\bar{X}} N^{n-1} \omega^n \, ,
\]
and we see that $N^{n-1} \omega^n$ is an invariant volume form. This proves
the following result.

{\sc Theorem 7.3.}(\cite{St}) {\it Condition (\ref{as2}) is sufficient for the 
existence of an invariant measure for the reduced Chaplygin equations (\ref{as}).}

{\sc Remark 7.4.} Stanchenko observes that if $F$ satisfies (\ref{as2}), 
the semi-basic character of both $\theta_{L^*}$ and $\Xi$ imply that $F$ is 
necessarily the pullback of a function on $Q/G (=M)$.

It turns out that condition (\ref{as2}) can be relaxed to some extent: 
it suffices to require the almost symplectic 2-form $\omega$ to be globally 
conformal symplectic, that is, that there exists a function 
$F \in C^{\infty}(T(Q/G))$ such that
\begin{equation}\label{as3}
dF \wedge \omega = - d\omega \, .
\end{equation}
Theorem 7.3 still holds in this case, with (\ref{as3}) replacing (\ref{as2}). 
The previous remark also remains valid: the function $F$ is necessarily the 
pullback of a function on $Q/G$. Note that (\ref{as2}) obviously implies (\ref{as3}). 

However, even the weaker condition (\ref{as3}) is not necessary in general
for the existence of an invariant volume form on $T(Q/G)$. 
To derive a necessary condition, let us suppose that $\mu$ is an invariant 
volume form for the dynamics $\bar{X}$ on $T(Q/G)$. We then necessarily have that 
$\mu = k \, \omega^n$, for some nowhere vanishing function $k \in C^{\infty}(T(Q/G))$.
Restricting ourselves to a connected component of $Q/G$ if
need be, we may always assume $k$ is strictly positive. It follows that
\begin{eqnarray*}
0 = {\cal L}_{\bar{X}} \mu &=& \bar{X}(k) \, \omega^n + k \, {\cal L}_{\bar{X}} \omega^n \\
&=& \bar{X}(k) \, \omega^n + n k \, {\cal L}_{\bar{X}} \omega \wedge \omega^{n-1} = \bar{X}(k) \, \omega^n - n k \, i_{\bar{X}} d \, \Xi \wedge \omega^{n-1} \, .
\end{eqnarray*}
The $2n$-form $i_{\bar{X}} d \, \Xi \wedge \omega^{n-1}$ determines a function 
$h \in C^{\infty}(T(Q/G))$ by
\begin{equation}\label{lah}
i_{\bar{X}} d \, \Xi \wedge \omega^{n-1} = \frac{h}{n} \, \omega^n \, .
\end{equation}
Therefore, we have
\begin{equation}\label{charac}
\bar{X}(k) = k h \quad \hbox{or, equivalently}, \quad \bar{X}(\hbox{ln}\,k) = h \, .
\end{equation}
This essentially yields the same characterization as the one derived in \cite{St}. 
Now, conversely, assume there exists a function $k$ satisfying (\ref{charac}),
with $h$ defined by (\ref{lah}). Going through the above computations in reverse
order then shows that the $2n$-form $k\,\omega^n$ is an invariant volume
form of $\bar{X}$. We may therefore conclude that the existence of a globally
defined function $k$ for which (\ref{charac}) holds is not only a necessary but 
also a sufficient condition for the existence of an invariant volume form. 
It is interesting to note that in \cite{St}, Stanchenko has proved that in case 
the reduced Lagrangian is of kinetic energy type, $L^* = \frac{1}{2} \tilde{g}$, 
and if there exists a solution $k$ of (\ref{charac}) which is basic, i.e.\ which 
is the pullback of a function on the base space $Q/G$, then the volume form 
$\mu = k \, \omega^n$ remains an invariant of the reduced dynamics if a potential 
energy function $\tilde{V} \in C^{\infty}(Q/G)$ is included in the Lagrangian $L^*$ 
(coming from a $G$-invariant potential added to the given Chaplygin system). 

Obviously, however, (\ref{charac}) is not a very handy criterium to deal with in 
practice. We will now see that, at least for a subclass of generalized Chaplygin 
systems, it may be replaced by a more easily manegeable condition.  

From (\ref{lah}), we can deduce a local expression for $h$. After some computations, 
we get
\[
h= \sum_{a,b} \tilde{g}^{ab} \frac{\partial \alpha_b}{\partial \dot{q}^a} \, ,
\]
and, using (\ref{local}), this further becomes
\[
h= \sum_{f,b}
\tilde{g}^{fb} \dot{q}^e \tilde{g}_{bc} (B^c_{ef} + B^c_{fe}) = 
\sum_{c,e} \dot{q}^e (B^c_{ec} + B^c_{ce}) \, .
\]
Note that $S^*dh$ is the pullback of a basic 1-form, i.e.\ $S^*dh = \tau_{Q/G}^* \beta$, 
where the local expression for $\beta$ reads 
\[
\beta = h_e(q)dq^e = \sum_{c} (B^c_{ec} + B^c_{ce})dq^e \, ,
\]
with $h_e = {\partial h}/{\partial \dot{q}^e}$. 
Let us assume now that there exists a basic function $k$ for which (\ref{charac}) 
holds. We then have that $S^*(d\bar{X}(\hbox{ln}\,k)) = d(\hbox{ln}\,k)$. Therefore, 
taking the differential of both hand-sides of (\ref{charac}) and applying $S^*$
to the resulting equality, we obtain
\[
d(\hbox{ln}\,k) = \beta \, ,
\]
or, in local coordinates,
\[
\frac{\partial (\hbox{ln}\,k)}{\partial q^e} = h_e(q) = \sum_{c} (B^c_{ec} + B^c_{ce}) \, , \; e=1,...,n \, .
\]
We thus see that, if (\ref{charac}) admits a solution $k$ which is basic, then 
the 1-form $\beta$ is exact. 

It turns out that, for systems for which $\tilde{V} = 0$
(i.e.\ $L^*$ is a pure kinetic energy Lagrangian), the previous result even has 
a more definitive character. Indeed, let $\tilde{V} = 0$ and suppose there exists a 
function $k \in C^{\infty}(T(Q/G))$ (not necessarily basic) satisfying (\ref{charac}). 
Then, we have that
\[
S^*d\bar{X}(\hbox{ln}\,k) = S^*dh \, .
\]
In local coordinates, this becomes
\[
\frac{\partial \bar{X}(\hbox{ln}\,k)}{\partial \dot{q}^e} = h_e(q) \, , \; e=1,...,n \, .
\]
But $\displaystyle{\bar{X}(\hbox{ln}\,k) = 
\dot{q}^a \frac{\partial (\hbox{ln}\,k)}{\partial q^a} 
+ \bar{X}^a \frac{\partial (\hbox{ln}\,k)}{\partial \dot{q}^a}}$, where $\bar{X}^a = 
-( \tilde{g}^{ab}\alpha_b + \dot{q}^b \dot{q}^c \tilde{\Gamma}_{bc}^a)$, and so we have
\[
\frac{\partial \bar{X}(\hbox{ln}\,k)}{\partial \dot{q}^e} = 
\frac{\partial (\hbox{ln}\,k)}{\partial q^e} + 
\dot{q}^a \frac{\partial^2 (\hbox{ln}\,k)}{\partial q^a \partial \dot{q}^e} + 
\frac{\partial \bar{X}^a}{\partial \dot{q}^e} \frac{\partial (\hbox{ln}\,k)}
{\partial \dot{q}^a} + 
\bar{X}^a \frac{\partial^2 (\hbox{ln}\,k)}{\partial \dot{q}^a \partial \dot{q}^e} \, .
\]
In points $0_q$ of the zero section of $T(Q/G)$ this reduces to
\[
\frac{\partial \bar{X}(\hbox{ln}\,k)}{\partial \dot{q}^e}\Big|_{0_q} =
\frac{\partial (\hbox{ln}\,k)}{\partial q^e}\Big|_{0_q}\,.
\]
If we now define the basic function $\tilde{k} = k \circ s$, where 
$s: Q/G \longrightarrow T(Q/G), q \longmapsto (q,0)$ determines the zero 
section, we derive from the above that
\[
\frac{\partial (\hbox{ln}\,\tilde{k})}{\partial q^e}(q) = 
\frac{\partial (\hbox{ln}\,k)}{\partial q^e}(q,0)= h_e(q) \, , \; e=1,...,n \, , 
\]
and, hence, it follows again that the 1-form $\beta$ is exact. 

Conversely, if the 1-form $\beta$ is exact, say $\beta = df$ for some function
$f \in C^{\infty}(Q/G)$, and putting $k = \exp({\tau_{Q/G}}^*f)$, one easily
verifies that $k$ satisfies (\ref{charac}). This obviously also holds in the presence 
of a potential (i.e.\ if $\tilde{V} \neq 0$).

Summarizing the above, we have proved the following interesting result.

{\sc Theorem 7.5.} {\it For a generalized Chaplygin system with a Lagrangian 
of kinetic energy type, i.e.\ $L=\frac{1}{2} g$, there exists an invariant volume 
form for the reduced dynamics $\bar{X}$ on $T(Q/G)$ iff the basic 1-form $\beta$, 
defined by $S^*dh = \tau^*_{Q/G}\beta$ (with $h$ given by (\ref{lah})), is exact. 
The `if' part also holds if $L$ is of the form $L = \frac{1}{2}g - V$,
with $V$ a $G$-invariant potential.}

Therefore, if we manage to find a particular example of a system with a kinetic 
energy type Lagrangian for which $\beta$ is not exact, we shall have proved that 
the answer to Koiller's question about the existence of an invariant volume form for
all generalized Chaplygin systems, is negative. In particular, for such a counter
example it suffices to show that the corresponding $\beta$ is not closed.

7.2. {\it A counter example}

Let us consider the following modified version of the classical example of 
the nonholonomic free particle \cite{BaSn}.

Consider a particle moving in space, so $Q = \hbox{\ddpp R}^3$,
subject to the non-holonomic constraint 
\[
\phi = \dot{z} - y x \dot{x} \, .
\]

The Lagrangian function is the kinetic energy corresponding to the standard metric 
$g = dx \otimes dx + dy \otimes dy + dz \otimes dz$. Therefore,  
\[
L = \frac{1}{2} \left(\dot{x}^2 + \dot{y}^2 + \dot{z}^2 \right) \; .
\]

The constraint submanifold is defined by the distribution
\[
{\cal D} = <\frac{\partial}{\partial x} + y x \frac{\partial}{\partial z},
\frac{\partial}{\partial y}> \, .
\]
  
Consider the Lie group $G=\hbox{\ddpp R}$ with its trivial 
action by translations on $Q$:
\[
\begin{array}{rccl}
\Phi : &  G \times Q & \longrightarrow & Q \\
& (s,(x,y,z)) & \longmapsto & (x,y,z+s) \, .
\end{array}
\]

Note that ${\cal D}$ is the horizontal subspace of a connection $\gamma$ 
on the principal fiber bundle $Q \longrightarrow Q/G$, where 
$\gamma = dz-yxdx $. Therefore, this system belongs to the class of generalized 
Chaplygin systems.

The curvature of $\gamma$ is given by
\[
\Omega^{\gamma} = x dx \wedge dy\, .
\]

The induced metric $\tilde{g}$ on $Q/G \cong \R^2$ is
\[
\tilde{g} = (1+x^2y^2)dx \otimes dx + dy \otimes dy \, .
\]

The metric connection tensor $\tilde{K}$ is determined by
\[
\tilde{K}(\frac{\partial}{\partial x},\frac{\partial}{\partial x},
\frac{\partial}{\partial y}) = x^2y \, , \; \;
\tilde{K}(\frac{\partial}{\partial y},\frac{\partial}{\partial x},
\frac{\partial}{\partial y}) = 0 \, .
\]

Then, the contorsion of the affine connection $\nabla^H_1$ here reads 
\[
B=x^2 y dx \otimes dx \otimes \frac{\partial}{\partial y} -
\frac{x^2 y}{1+x^2y^2} dx \otimes dy \otimes \frac{\partial}{\partial x} \, .
\]

Finally, the 1-form $\beta$ associated to the reduced Chaplygin system is given by
\[
\beta = - \frac{x^2 y}{1+x^2 y^2} dy \, ,
\]
which is clearly not closed. Hence, according to Theorem 7.5, there is no invariant
volume form for the reduced dynamics. 

Note that in this example the distribution $\cal D$ has `length' $1$ at all points
of $\hbox{\ddpp R}^3$ not belonging to the plane $x=0$, since 
${\cal D}+[{\cal D},{\cal D}]$ spans the full tangent space at all points for 
which $x \neq 0$. We leave it as a challenge for the reader to find a more generic 
example (on $\hbox{\ddpp R}^3$) for which the constraint distribution has length $1$ 
everywhere.

{\sc Remark 7.6.} The classical model for the nonholonomic free particle corresponds 
to the constraint $\phi=\dot{z}-y\dot{x}$. Here the constraint distribution does have
length $1$ everywhere, however, after performing the appropriate computations, 
we obtain $\displaystyle{\beta=-\frac{y}{1+y^2}dy}$, which is clearly exact: 
$\beta = df$, with $\displaystyle{f = \frac{1}{2} \hbox{ln} \left( \frac{1}{1+y^2} \right)}$. Then, $\displaystyle{k=\frac{1}{\sqrt{1+y^2}}}$ and the invariant measure defined by 
$k\, \omega^n$ leads, using Euler's integrating factor technique, to the constant of 
the motion
\[
\varphi= \dot{x} \sqrt{1+y^2} \, 
\]
which was also obtained by different methods in \cite{BaGrDo,BlKrMaMu,CoLe}.

\begin{center}
{\it Appendix: Metric connections}
\end{center}
\setcounter{equation}{0}
In this section, we want to collect some simple facts about metric
connections that have been useful for the formulation of generalized
Chaplygin systems. 

Let $Q$ be a $n$-dimensional differentiable manifold, the configuration
space with Riemannian metric ${g}$. We denote by $\nabla^{g}$ the
Levi-Civita or Riemannian connection associated to the metric ${g}$. 

{\sc Definition A.1.} {\it An affine connection $\nabla$ is called metric with 
respect to $g$ if $\nabla g=0$, that is, 
\[
Z(g(X,Y)) = g(\nabla_Z X,Y)+g(X,\nabla_Z Y) \, ,
\]
for all $X,Y,Z \in {\frak X}(Q)$.}

Let $\nabla$ be a metric connection with respect to $g$. In the
following proposition, we prove that $\nabla$ is determined by its
torsion $T$. 

{\sc Proposition A.2.} {\it Let $T$ be a skewsymmetric (1,2) tensor on $Q$. 
Then there exists a unique metric connection $\nabla$ whose torsion is precisely $T$.}

{\it Proof.} Let us suppose that there exists such metric connection $\nabla$. Then
we have that 
\begin{eqnarray*}
Z(g(X,Y)) &=& g(\nabla_Z X,Y)+g(X,\nabla_Z Y) \, , \\
X(g(Z,Y)) &=& g(\nabla_X Z,Y)+g(Z,\nabla_X Y) \, , \\
Y(g(X,Z)) &=& g(\nabla_Y X,Z)+g(X,\nabla_Y Z) \, ,
\end{eqnarray*}
for all $X,Y,Z \in {\frak X}(Q)$. Now
\begin{eqnarray*}
&& Z(g(X,Y)) + X(g(Z,Y)) - Y(g(X,Z)) \\
&=& g(\nabla_X Z+\nabla_Z X,Y) + g(\nabla_Z Y-\nabla_Y Z,X) +
g(\nabla_X Y-\nabla_Y X,Z)\\ 
&=& g(2\nabla_X Z+T(Z,X)+[Z,X],Y) + g(T(Z,Y)+[Z,Y],X) + g(T(X,Y)+[X,Y],Z) \\
&=& 2g(\nabla_X Z,Y) + g(T(Z,X)+[Z,X],Y) + g(T(Z,Y)+[Z,Y],X) +
g(T(X,Y)+[X,Y],Z) \, . 
\end{eqnarray*}
Consequently, the connection $\nabla$ is uniquely determined by the
formula 
\begin{equation}
g(\nabla_X Z,Y)=g(\nabla^g_X Z,Y)-\frac{1}{2} \left(
g(Y,T(X,Z))+g(X,T(Z,Y))+g(Z,T(X,Y)) \right) \, . 
\end{equation}
\QED

This proposition implies that the Christoffel symbols
$\bar{\Gamma}^A_{BC}$ of the metric connection $\nabla$ in a local
chart $(q^A)$ are given by 
\[
\bar{\Gamma}^A_{BC} = \Gamma^A_{BC} - \frac{1}{2} g^{AK} \left( g_{KM}
T^M_{BC} + g_{BM} T^M_{CK} + g_{CM} T^M_{BK} \right) \, ,
\]
where $\Gamma^A_{BC}$ are the Christoffel symbols of the connection
$\nabla^g$ and $T=T^C_{AB} dq^A \otimes dq^B \otimes \frac{\partial }{\partial q^C}$. 

Another way to characterize metric connections is the following. Given
any affine connection on $Q$, we know that 
\[
\nabla_X Y = \nabla^g_X Y +S(X,Y) \, ,
\]
where $S$ is a (1,2) tensor field, called the contorsion. If $\nabla$ is a metric 
connection, then 
\begin{eqnarray*}
Z(g(X,Y)) &=& g(\nabla_Z X,Y)+g(X,\nabla_Z Y) \\
&=& g(\nabla^g_Z X + S(Z,X),Y)+g(X,\nabla^g_Z Y+S(Z,Y)) \\
&=& Z(g(X,Y)) + g(S(Z,X),Y) + g(X,S(Z,Y)) \, ,
\end{eqnarray*}
which implies that $g(S(Z,X),Y) + g(X,S(Z,Y))=0$. Then we have proved
the following 

{\sc Proposition A.3.} $\nabla$ {\it is a metric connection if and only if}
\begin{equation}\label{*}
g(S(Z,X),X)=0 \; \; \forall X,Z \in {\frak X}(Q) \, .
\end{equation}

As a consequence of the two characterizations we have obtained for
metric connections, we can establish the next result. 

{\sc Corollary A.4.} {\it There is a one-to-one correspondence between (1,2) 
tensors $S$ verifying (\ref{*}) and skewsymmetric (1,2) tensors $T$. This
correspondence is given by 
\[
S \longrightarrow T \, ,
\]
where $T(X,Y)=S(X,Y)- S(Y,X)$ and
\[
T \longrightarrow S \, ,
\]
where $g(S(X,Z),Y)=
-\displaystyle{\frac{1}{2}} \left((g(Y,T(X,Z))+g(X,T(Y,Z))+g(Z,T(X,Y)) \right)$.}

The equations for the geodesics of a metric connection are
\[
\nabla_{\dot{c}(t)} \dot{c}(t) = 0 \Longleftrightarrow
\nabla^g_{\dot{c}(t)} \dot{c}(t) = -S(\dot{c}(t),\dot{c}(t)) \, , 
\]
or, in local coordinates
\[
\ddot{q}^A + \Gamma^A_{BC} \dot{q}^B \dot{q}^C = \sum_{B<C} g^{AK} \left( 
g_{BM}T^M_{CK} + g_{CM}T^M_{BK} \right) \dot{q}^B \dot{q}^C \, ,
\]
for each $A=1,...,n$.

Finally, it is very important to note that the metric connections
preserve obviously the kinetic energy $K$ of the metric $g$, that is,
if $c(t)$ is a geodesic of $\nabla$, we have that 
\[
\frac{d}{dt}\left( \frac{1}{2} g(\dot{c}(t),\dot{c(t}) \right) =
g(\nabla_{\dot{c}(t)}\dot{c}(t),\dot{c}(t)) = 0 \, . 
\]

{\sc Acknowledgements.}
This work was partially supported by grant DGICYT (Spain) PB97-1257. 
J. Cort\'es wishes to thank the Spanish Ministerio de 
Educaci\'on y Cultura for a FPU grant and both J. Cort\'es and D. Mart\'{\i}n
de Diego wish to thank the Department of Mathematical Physics and Astronomy of 
the University of Ghent for its kind hospitality. We would like to thank 
S. Mart{\'\i}nez for helpful conversations and the referee for several useful comments.

\end{document}